\documentclass[11pt]{article}\textwidth 160mm\textheight 235mm
\usepackage{amsfonts}
\usepackage{amssymb}
\usepackage{graphicx}
\usepackage{amsmath}
\usepackage{tikz}
\usepackage{cancel}
\usepackage{todonotes}
\usetikzlibrary{matrix}
\usetikzlibrary{arrows,shapes}

\def\ve{\varepsilon}
\def\epsilon{\varepsilon}

\def\hat{\widehat}
\def\tilde{\widetilde}
\def\emp{\emptyset}

\def\dist{{\rm dist}}
\def\dom{{\rm dom}\,}
\def\span{{\rm span}\,}
\def\epi{{\rm epi\,}}
\def\rge{{\rm rge\,}}

\def\p{^{\prime}}

\def\N{{\cal N}}

\def\Lm{{\Lambda}}
\def\O{{\cal O}}

\def\d{{\rm d}}
\def\sub{\partial}
\def\B{\mathbb B}

\def\ox{\overline{x}}
\def\oy{\overline{y}}
\def\oz{\overline{z}}

\def\disp{\displaystyle}

\def\tto{\;{\lower 1pt\hbox{$\rightarrow$}}\kern-10pt
\hbox{\raise 2pt\hbox{$\rightarrow$}}\;}
\def\Hat{\widehat}
\def\Tilde{\widetilde}
\def\Bar{\overline}
\def\ra{\rangle}
\def\la{\langle}
\def\ve{\varepsilon}
\def\B{I\!\!B}
\def\h{\hfill\Box}
\def\R{\mathbb{R}}
\def\N{I\!\!N}
\def\ox{\bar{x}}
\def\oy{\bar{y}}
\def\oz{\bar{z}}
\def\ov{\bar{v}}

\def\ou{\bar{u}}

\def\op{\bar{p}}

\def\co{\mbox{\rm co}\,}
\def\cone{\mbox{\rm cone}\,}

\def\int{\mbox{\rm int}\,}
\def\gph{\mbox{\rm gph}\,}
\def\epi{\mbox{\rm epi}\,}

\def\dom{\mbox{\rm dom}\,}

\def\ker{\mbox{\rm ker}\,}

\def\aff{\mbox{\rm aff}\,}

\def\h{\hfill\triangle}
\def\dn{\downarrow}
\def\O{\Omega}
\def\ph{\varphi}
\def\emp{\emptyset}
\def\st{\stackrel}
\def\oR{\Bar{\R}}
\def\lm{\lambda}
\def\gg{\gamma}
\def\dd{\delta}
\def\al{\alpha}
\def\Th{\Theta}
\def\N{I\!\!N}
\def\th{\theta}

\def\sce{\setcounter{equation}{0}}
\begin{document}
\begin{center}
{\bf CRITICAL MULTIPLIERS IN VARIATIONAL SYSTEMS\\ VIA SECOND-ORDER GENERALIZED DIFFERENTIATION}\\[2ex]
BORIS S. MORDUKHOVICH\footnote{Department of Mathematics, Wayne State University, Detroit, MI 48202, USA and the RUDN University, Moscow 117198, Russia (boris@math.wayne.edu). Research of this author was partly supported by the National Science Foundation under grants DMS-1007132 and DMS-1512846, by the Air Force Office of Scientific Research under grant \#15RT0462, and by the Ministry of Education and Science of the Russian Federation (Agreement number 02.a03.21.0008 of 24 June 2016).} and M. EBRAHIM SARABI\footnote{Department of Mathematics, Miami University, Oxford, OH 45065, USA (sarabim@miamioh.edu).}
\end{center}
\small{\bf Abstract.} In this paper we introduce the notions of critical and noncritical multipliers for variational systems and extend to a general framework the corresponding notions by Izmailov and Solodov developed for classical Karush-Kuhn-Tucker (KKT) systems. It has been well recognized that critical multipliers are largely responsible for slow convergence of major primal-dual algorithms of optimization. The approach of this paper allows us to cover KKT systems arising in various classes of smooth and nonsmooth problems of constrained optimization including composite optimization, minimax problems, etc. Concentrating on a polyhedral subdifferential case and employing recent results of second-order subdifferential theory, we obtain complete characterizations of critical and noncritical multipliers via the problem data. It is shown that noncriticality is equivalent to a certain calmness property of a perturbed variational system and that critical multipliers can be ruled out by full stability of local minimizers in problems of composite optimization. For the latter class we establish the equivalence between noncriticality of multipliers and robust isolated calmness of the associated solution map and then derive explicit characterizations of these notions via appropriate second-order sufficient conditions. It is finally proved that the Lipschitz-like/Aubin property of solution maps yields their robust isolated calmness.\\[1ex]
{\bf Keywords} Variational systems, Composite optimization, Critical and noncritical multipliers, Generalized differentiation, Piecewise linear functions, Robust isolated calmness, Lipschitzian stability\\[1ex]
{\bf Mathematics Subject Classification (2000)} 90C31, 49J52, 49J53

\newtheorem{Theorem}{Theorem}[section]
\newtheorem{Proposition}[Theorem]{Proposition}
\newtheorem{Remark}[Theorem]{Remark}
\newtheorem{Lemma}[Theorem]{Lemma}
\newtheorem{Corollary}[Theorem]{Corollary}
\newtheorem{Definition}[Theorem]{Definition}
\newtheorem{Example}[Theorem]{Example}
\renewcommand{\theequation}{{\thesection}.\arabic{equation}}
\renewcommand{\thefootnote}{\fnsymbol{footnote}}

\normalsize\vspace*{-0.15in}
\section{Introduction}\sce\vspace*{-0.05in}

In recent years it has been well recognized that the so-called {\em critical Lagrange multipliers} play a serious {\em negative} role in the convergence of primal-dual algorithms of numerical optimization. Namely, their existence implies a slow (less than superlinear) convergence of major algorithms of such types. We refer the reader to the monograph by Izmailov and Solodov \cite{is14} and their excellent survey paper \cite{is15} devoted to comprehensive discussions of this phenomenon in problems of nonlinear programming (NLPs) with ${\cal C}^2$-smooth data; see also the experts' comments to \cite{is15} and the authors' rejoinder published in the same issue of TOP.

The main interest of this paper is to introduce, characterize, and apply critical multipliers and their noncritical counterparts for the following class of {\em variational systems} of the {\em subdifferential type}. Given mappings $\Phi\colon\R^n\to\R^m$, $\Psi\colon\R^n\times\R^m\to\R^l$ and an extended-real-valued function $\th\colon\R^m\to\oR:=(-\infty,\infty]$, consider the system of equations and inclusions defined by
\begin{equation}\label{VS}
\Psi(x,v)=0,\;v\in\partial\th\big(\Phi(x)\big),
\end{equation}
where $\partial\th$ stands for an appropriate subdifferential of $\th$. In this paper we mainly deal with {\em convex} functions $\th$, and so their subdifferential is in the classical sense of convex analysis.

Note that \eqref{VS}, being applied to optimization problems, can be treated as a ``generalized KKT system." Indeed, consider the following problem of {\em composite optimization}:
\begin{equation}\label{comp}
\mbox{minimize }\;\ph(x):=\ph_0(x)+\theta\big(\Phi(x)\big),\quad x\in\R^n.
\end{equation}
Although \eqref{comp} is written in the unconstrained format, it implicitly includes the constraints $\Phi(x)\in\dom\th:=\{z\in\R^m|\;\th(z)<\infty\}$. This model has been widely used as a convenient form to study of various classes of constrained optimization problems, which go far beyond usual inequality and equality constraints in nonlinear programming; see, e.g., \cite{mrs,rw} for more discussions and references. Denoting the {\em Lagrangian} of \eqref{comp} by
\begin{equation}\label{lagr}
L(x,v):=\ph_0(x)+\la\Phi(x),v\ra,\quad(x,v)\in\R^n\times\R^m,
\end{equation}
and choosing $\Psi:=\nabla_x L$, we see that \eqref{VS} reduces to the KKT system for the composite optimization problem \eqref{comp} and thus for its more conventional specifications.

In this paper we define the notions of critical and noncritical multipliers for variational systems \eqref{VS} and conduct a rather comprehensive study of them in the case where the mappings $\Phi$ and $\Psi$ are sufficiently smooth while the convex function $\th$ is generally extended-real-valued (hence definitely nonsmooth) but piecewise linear, i.e., its epigraph is a convex polyhedron. We indicate that $\th$ belongs to this class of {\em convex piecewise linear} functions by writing $\th\in CPWL$.

Our analysis and applications of criticality and noncriticality in the framework of \eqref{VS} with $\th\in CPWL$ are heavily based on the recently developed {\em second-order subdifferential} calculations for this class of functions \cite{ms15}, which allow us to efficiently characterize such multipliers entirely in terms of the given polyhedron data and then constructively apply them to the study of some important notions of stability for variational systems and optimization problems.

As mentioned above, critical multipliers have a negative influence on the convergence rate for major primal-dual algorithms of optimization. Thus it is crucial  from the computational viewpoint to recognize situations where critical multipliers {\em cannot} be associated with a particular local minimizer and then to develop algorithms which perform well in searching not arbitrary but such ``good" optimal solutions. It has been conjectured by the first author \cite{m15} that the property of {\em full stability} \cite{lpr} of a local minimizer rules out the existence of critical multipliers associated with this minimizer. It has also been conjectured in \cite{m15} that even the weaker property of {\em tilt stability} \cite{pr98} would exclude the existence of critical multipliers under appropriate assumptions. Some results on these conjectures for NLPs are obtained in \cite{gm15,iz15,mn15} and are discussed in the sequel together with new developments concerning the composite optimization model \eqref{comp}. Note that resolving these conjectures in the affirmative opens the gate to constructively verify the possibility of ruling out critical multipliers in practical situations, since by now we have efficient second-order characterizations of full and tilt stability for large classes of optimization and variational problems; see more discussions and references below.

Another benefit of the obtained characterizations of noncriticality, which has never been exploited before in the literature, is establishing the equivalence between noncriticality of a Lagrange multiplier in composite optimization and {\em robust isolated calmness} of the solution map to the corresponding {\em canonically perturbed} KKT system. The latter term has been recently coined in \cite{dsz} to distinguish this robust notion from (nonrobust) isolated calmness and the equivalent strong metric subregularity of the inverse; see, e.g., \cite{dr}. In this way we derive a new {\em second-order characterization} of robust isolated calmness for the KKT system associated with \eqref{comp} while expressing it entirely in terms of the given data. It is shown finally that the robust isolated calmness of the latter system is implied by its Lipschitz-like/Aubin property.

The rest of the paper is organized as follows. Section~2 briefly recalls some tools and results of variational analysis and generalized differentiation widely used in the subsequent material. In Section~3 we define critical and noncritical multipliers for \eqref{VS}, establish equivalent descriptions of {\em critical multipliers} for the case of $\th\in CPWL$, and specify them for particular KKT systems in smooth and nonsmooth optimization. Section~4 is mainly devoted to characterizing {\em noncritical multipliers} for the variational system \eqref{VS} with $\th\in CPWL$ via a certain calmness property (defined in this paper as ``semi-isolated calmness") that involves the solution map to a canonically perturbed counterpart of \eqref{VS}. We also present here a new {\em second-order sufficient condition} (SOSC) for noncriticality in the general framework of \eqref{VS} and show that it ensures the {\em strict local optimality} in composite optimization problems.

Starting with Section~5, we focus solely on the composite model of optimization \eqref{comp} with $\th\in CPWL$ therein. Section~5 justifies for this model the aforementioned conjecture on {\em ruling out} the existence of critical minimizers associated with {\em fully stable} local minimizers of \eqref{comp}. In Section~6 we discuss some qualification conditions allowing us to exclude critical minimizers associated with {\em tilt-stable} multipliers while, on the other hand, present examples showing that generally it is not the case in various settings of NLP.

Section~7 is devoted to the study of {\em isolated calmness} and its {\em robust} counterpart for the solution map to the canonically perturbation of \eqref{VS} and its KKT specification for \eqref{comp}. By implementing a new approach based on the developed critical multiplier theory, we establish close relationships between noncriticality and isolated calmness for general systems \eqref{VS} and then strengthen them for the case of $\th\in CPWL$ with applications to composite optimization. This approach allows us, in particular, to {\em characterize} both isolated calmness and its robust version for the KKT system associated with a locally optimal solution to \eqref{comp} when $\th\in CPWL$  by the corresponding specification of the {\em SOSC} for noncritical multipliers established in Section~4.

Section~8 justifies the validity of the noncriticality, nondegeneracy, and robust isolated calmness properties of KKT solution maps in composite optimization under their Lipschitz-like stability. The concluding Section~9 contains final discussions with the emphasis on the major points of the paper and formulations of some open questions of the future research.

Throughout the paper we use the standard notation from variational analysis; cf.\ \cite{m06,rw}. Recall that $\B_r(x)$ stands for the closed ball
centered at $x$ with radius $r>0$, while $\B$ indicates the closed unit ball in the space in question if no confusion arises.\vspace*{-0.15in}

\section{Preliminaries from Variational Analysis}\sce\vspace*{-0.05in}

In this section we first briefly review, following mainly the books \cite{m06,rw}, basic constructions of variational analysis and generalized differentiation employed in the paper and then recall some recent results of \cite{ms15} concerning CPWL functions that are largely used in what follows.

Given a set $\O\subset\R^m$, its (Fr\'echet) {\em regular normal cone} is defined by
\begin{equation}\label{2.1}
\Hat N(z;\O):=\disp\Big\{v\in\R^m\Big|\;\limsup_{u\st{\O}{\to}z}\frac{\la v,u-z\ra}{\|u-z\|}\le 0\Big\},\quad z\in\O,
\end{equation}
where the symbol $u\st{\O}{\to}z$ means that $u\to z$ with $u\in\O$. Construction \eqref{2.1} is also called the ``prenormal cone" to $\O$ at $z$ due to the fact that it fails to possess some expected properties of normals to closed sets being often empty at boundary points as, e.g., for $\O:=\epi(-|x|)\subset\R^2$ at $z=(0,0)$. The (Mordukhovich) {\em limiting normal cone} to $\O$ at $\oz$ defined by
\begin{equation}\label{2.2}
N(\oz;\O)=\big\{v\in\R^m\big|\;\exists\,z_k\st{\O}{\to}\oz,\;v_k\in\Hat N(z_k;\O)\;\mbox{ with }\;v_k\to v\;\mbox{ as }\;k\to\infty\big\}
\end{equation}
possesses the aforementioned and other required properties of generalized normals and, despite its nonconvexity, enjoys---together with the associated subdifferential and coderivative constructions for extended-real-valued functions and set-valued mappings/multifunctions, respectively,---comprehensive calculus rules based on variational/extremal principles of variational analysis. If $\O$ is convex, both constructions \eqref{2.1} and \eqref{2.2} reduce to the classical normal cone of convex analysis. Recall the duality relationship
\begin{eqnarray*}
\Hat N(z;\O)=T(z;\O)^*:=\big\{v\in\R^m\big|\;\la v,w\ra\le 0\;\mbox{ for all }\;w\in T(z;\O)\big\}
\end{eqnarray*}
between \eqref{2.1} and the (Bouligand-Severi) {\em tangent cone} $T(z;\O)$ to $\O$ at $z\in\O$ defined by
\begin{eqnarray}\label{2.5}
T(z;\O):=\big\{w\in\R^m\big|\;\exists\,z_k\st{\O}{\to}z,\;\al_k\ge 0\;\mbox{ with }\;\al_k(z_k-z)\to w\;\mbox{ as }\;k\to\infty\big\}.
\end{eqnarray}

For an extended-real-valued function $\th\colon\R^m\to\oR$, consider the two limiting subdifferential constructions associated with \eqref{2.2}: the {\em basic subdifferential} and the {\em singular subdifferential} of $\th$ at $\oz\in\dom\th$ given, respectively, by
\begin{equation}\label{2.6}
\partial\th(\oz):=\big\{v\in\R^m\big|\;(v,-1)\in N\big((\oz,\th(\oz));\epi\th)\big)\big\},
\end{equation}
\begin{equation}\label{2.66}
\partial^\infty\th(\oz):=\big\{v\in\R^m\big|\;(v,0)\in N\big((\oz,\th(\oz));\epi\th)\big)\big\}.
\end{equation}
We know that for convex functions $\th$ the basic subdifferential (\ref{2.6}) agrees with the subdifferential of convex analysis and that for the general class of lower semicontinuous (l.s.c.) functions $\th$  the singular subdifferential \eqref{2.66} reduces to $\{0\}$ if and only if $\th$ is locally Lipschitzian around $\oz$. Note also that we have the representation $\partial^\infty\th(\oz)=N(\oz;\dom\th)$ for convex functions $\th$ and that
\begin{equation*}
N(\oz;\O)=\sub\dd(\oz;\O)=\sub^\infty\dd(\oz;\O),\quad\oz\in\O,
\end{equation*}
for any set $\O$ via its indicator function $\dd_\O(z)=\dd(z;\O):=0$ for $z\in\O$ and $\dd(z;\O):=\infty$ otherwise.

Consider next a set-valued mapping $F\colon\R^n\tto\R^p$ with its domain and graph given by
$$
\dom F:=\big\{x\in\R^n\big|\;F(x)\ne\emp\big\}\quad\mbox{and}\quad\gph F:=\big\{(x,y)\in\R^n\times\R^p\big|\;x\in F(x)\big\}
$$
and define for it the following generalized differential notions via tangential and normal constructions from \eqref{2.1}--\eqref{2.5} to its graph. The {\em regular coderivative} and the {\em limiting coderivative} to $F$ at $(\ox,\oy)\in\gph F$ are given, respectively, by
\begin{equation}\label{2.7}
\Hat D^*F(\ox,\oy)(v):=\big\{u\in\R^n\big|\;(u,-v)\in\Hat N\big((\ox,\oy);\gph F\big)\big\},\quad v\in\R^p,
\end{equation}
\begin{equation}\label{2.8}
D^*F(\ox,\oy)(v):=\big\{u\in\R^n\big|\;(u,-v)\in N\big((\ox,\oy);\gph F\big)\big\},\quad v\in\R^p,
\end{equation}
while the {\em graphical derivative} of $F$ at $(\ox,\oy)$ is defined by
\begin{equation}\label{gder}
DF(\ox,\oy)(u):=\big\{v\in\R^p\big|\;(u,v)\in T\big((\ox,\oy);\gph F\big)\big\},\quad u\in\R^n.
\end{equation}
If $F\colon\R^n\to\R^p$ is single-valued, we drop $\oy$ in the notation \eqref{2.7}--\eqref{gder}. The smoothness of $F$ around $\ox$ in the latter case yields the representations
\begin{eqnarray*}
DF(\ox)(u)=\big\{\nabla F(\ox)u\big\},\quad\hat D^*F(\ox)(v)=D^*F(\ox)(v)=\big\{\nabla F(\ox)^*v\big\}\;\mbox{ for }\;u\in\R^n,\;v\in\R^p,
\end{eqnarray*}
where the symbol $A^*$ for the matrix $A$ signifies the matrix transposition/adjont operator.

In what follows we often use the mappings $D\partial\th$ and $D^*\partial\th$, which are constructions of {\em second-order} generalized differentiation for extended-real-valued functions $\th\colon\R^m\to\oR$ via the ``derivative-of-derivative" approach developed in \cite{m06} for the case of coderivatives. Note that  for functions $\th$ of class ${\cal C}^2$ near $\oz$ it holds
\begin{equation*}
\big(D\partial\th\big)\big(\oz,\th(\oz)\big)(u)=\big(D^*\partial\th\big)\big(\oz,\th(\oz)\big)(u)=\big\{\nabla^2\th(\oz)u\big\},\quad u\in\R^m,
\end{equation*}
due to the classical Hessian symmetry, while it is not the case for more general functions and also in infinite dimensions. Efficient applications of the aforementioned second-order constructions given in this paper for the case of $\th\in CPWL$ are largely based on the obtained second-order calculations of these constrictions for such functions $\th$ entirely in terms of their initial data.\vspace*{0.05in}

Dealing with {\em convex piecewise linear} functions $\th\colon\R^m\to\oR$, $\th\in CPWL$, recall their following equivalent descriptions taken from \cite[Theorem~2.49]{rw}:\\[1ex]
$\bullet$ The epigraphical set $\epi\th$ is a convex polyhedron in $\R^{m+1}$.\\[1ex]
$\bullet$ There are $\alpha_i\in\R$, $l\in\N$, and $a_i\in\R^m$ for $i\in T_1\colon=\{1,\ldots,l\}$ such that $\th$ is represented by
\begin{equation}\label{eq00}
\th(z)=\max\big\{\la a_1,z\ra-\alpha_1,\ldots,\la a_l,z\ra-\alpha_l\big\}\;\mbox{ if }\;z\in\dom\th
\end{equation}
and $\th(z)=\infty$ otherwise, where the domain set $\dom\th$ is a convex polyhedron given by
\begin{equation}\label{dom}
\dom\th=\big\{z\in\R^m\big|\;\la d_i,z\ra\le\beta_i\;\mbox{ for all }\;i\in T_2:=\{1,\ldots,p\}\big\}
\end{equation}
with some elements $d_i\in\R^m$, $\beta_i\in\R$, and $p\in\N$.

It follows from (\ref{eq00}) that each $\th\in CPWL$ can be expressed in the summation form
\begin{equation}\label{cpwl1}
\th(z)=\max\big\{\la a_1,z\ra-\alpha_1,\ldots,\la a_l,z\ra-\alpha_l\big\}+\dd(z;\dom\th),\quad z\in\R^m.
\end{equation}
It is observed in \cite[Proposition~3.2]{ms15} that, besides (\ref{dom}), the domain of $\th$ admits the representation
$\dom\th=\bigcup^{l}_{i=1}{C_i}$ with $l$ taken from (\ref{eq00}) and the sets $C_i$, $i\in T_1$, defined by
\begin{equation}\label{pwlr1}
C_i:=\big\{z\in\dom\th\big|\;\la a_j,z\ra-\al_j\le\la a_i,z\ra-\al_i,\;\;\mbox{for all}\;\;j\in T_1\big\}.
\end{equation}
Consider now the corresponding active index subsets in (\ref{pwlr1}) and (\ref{dom}) given by
\begin{equation}\label{active2}
K(\oz):=\big\{i\in T_1\big|\;\oz\in C_i\big\}\;\mbox{ and }\;I(\oz):=\big\{i\in T_2\big|\;\la d_i,\oz\ra=\beta_i\big\}
\end{equation}
and recall the formula for $\partial\th(\oz)$ at $\oz\in\dom\th$ obtained in \cite[Proposition~3.3]{ms15}:
\begin{equation}\label{fos}
\partial\th(\oz)=\co\big\{a_i\big|\;i\in K(\oz)\big\}+N(\oz;\dom\th)=\co\big\{a_i\big|\;i\in K(\oz)\big\}+\cone\big\{d_i\big|\;\;i\in I(\oz)\big\},
\end{equation}
where ``co" and ``cone" stand for the convex and conic hulls, respectively. Then for any $(\oz,\ov)\in\gph\partial\th$ we get from \eqref{fos} that $\ov=\ov_1+\ov_2$ with
\begin{equation}\label{eq06}
\disp\ov_1=\sum_{i\in K(\oz)}\bar\lm_i a_i\;\mbox{ with }\;\disp\sum_{i\in K(\oz)}\bar\lm_i=1,\;\bar\lm_i\ge 0\;\mbox{ and }\;
\disp\ov_2=\sum_{i\in I(\oz)}\bar\mu_id_i\;\mbox{ with }\;\bar\mu_i\ge 0.
\end{equation}
Recall also the well-known tangent cone representation
\begin{equation}\label{tanc}
T(\oz;\dom\th)=\big\{z\in\R^m\big|\;\la d_i,z\ra\le 0\;\;\mbox{for all}\;\;i\in I(\oz)\big\},\quad\oz\in\dom\th.
\end{equation}
Corresponding to (\ref{eq06}), define further the index subsets of positive multipliers for the given vectors $\ov_1$ and $\ov_2$ from (\ref{eq06}) by
\begin{equation}\label{eq05}
J_+(\oz,\ov_1):=\big\{i\in K(\oz)\big|\;\bar\lm_i>0\big\},\quad J_+(\oz,\ov_2):=\big\{i\in I(\oz)\big|\;\bar\mu_i>0\big\}
\end{equation}
and then consider the following sets constructed entirely in terms of the parameters in \eqref{eq00} and \eqref{dom} along arbitrary index subsets $P_1\subset Q_1\subset T_1$, $P_2\subset Q_2\subset T_2$ by\vspace*{-0.05in}
\begin{eqnarray}\label{eq080}
\begin{array}{lll}
{\cal F}_{\tiny\{P_1,Q_1\},\{P_2,Q_2\}}:&=\span\big\{a_i-a_j\big|\;i,j\in P_1\big\}+\cone\big\{a_i-a_j\big|\;(i,j)\in(Q_1\setminus P_1)\times P_1\big\}\\
&+\span\big\{d_i\big|\;i\in P_2\big\}+\cone\big\{d_i\big|\;i\in Q_2\setminus P_2\big\},
\end{array}
\end{eqnarray}\vspace*{-0.25in}
\begin{eqnarray}\label{eq081}
\begin{array}{ll}
{\cal G}_{\tiny\{P_1,Q_1\},\{P_2,Q_2\}}:=\Big\{u\in\R^n\Big|&\la a_i-a_j,u\ra=0\;\mbox{ if }\;i,j\in P_1,\\
&\la a_i-a_j,u\ra\le 0\;\mbox{ if }\;(i,j)\in(Q_1\setminus P_1)\times P_1,\\
&\la d_i,u\ra=0\;\mbox{ if }\;i\in P_2,\;\mbox{ and }\;\la d_i,u\ra\le 0\;\mbox{ if }\;i\in Q_2\setminus P_2\;\Big\}.
\end{array}
\end{eqnarray}\vspace*{-0.1in}
It is easy to deduce from the classical Farkas Lemma that
\begin{equation}\label{eq082}
{\cal G}_{\tiny\{P_1,Q_1\},\{P_2,Q_2\}}^*={\cal F}_{\tiny\{P_1,Q_1\},\{P_2,Q_2\}}\;\mbox{ for any }\;P_1\subset Q_1\subset T_1\;\mbox{ and }\;P_2\subset Q_2\subset T_2.
\end{equation}\vspace*{-0.15in}

We finish this section with the following extension of \cite[Theorem~3.4]{ms15} needed in the sequel. The additional information presented below can be deduced from the proof given therein. It is worth mentioning here that although the theorem presented below as well as Proposition~\ref{ccid} are written for the selected representation of $\ov$ in (\ref{eq06}), they are {\em invariant} with respect to different choices of $\ov_1$, $\ov_2$, $\bar\lm_i$, and $\bar\mu_i$ in (\ref{eq06}). Indeed, it has been shown in \cite[Propsoition~4.4]{ms15} that the second-order constructions are invariant with respect to such different choices. Furthermore, a close look at the proof of \cite[Theorem~3.4]{ms15} reveals that the obtained neighborhood $O$ in the following theorem is not depending on representation (\ref{eq06}). \vspace*{-0.05in}

\begin{Theorem}{\bf(description of points in the subdifferential graph of CPWL functions).}\label{gphpar} Let $\th\in CPWL$ with $(\oz,\ov)\in\gph\partial\th$. Then there is a neighborhood $O$ of $(\oz,\ov)$ such that for any $(z,v)\in(\gph\partial\th)\cap O$ we have $J_+(\oz,\ov_1)\subset K(z)$ and $J_+(\oz,\ov_2)\subset I(z)$, where $\ov_1$ and $\ov_2$  are taken from {\rm(\ref{eq06})}, and where $J_+(\oz,\ov_1)$ and $J_+(\oz,\ov_2)$ are defined in {\rm(\ref{eq05})}.
\end{Theorem}\vspace*{-0.3in}

\section{Critical Multipliers: Definition, Descriptions and Examples}\sce\vspace*{-0.05in}

In this section we define critical (and noncritical) multipliers for variational systems of type \eqref{VS}, establish its equivalent descriptions for the major case of $\th\in CPWL$ of our study and applications in this paper, and present several examples of multiplier criticality for particular classes of smooth and nonsmooth optimization problems.

In further developments and applications we impose the following connection between the mappings $\Phi$ and $\Psi$ in \eqref{VS} formulated via a given mapping $f\colon\R^n\to\R^n$. Assuming that $\Phi\colon\R^n\to\R^m$ is smooth, define $\Psi\colon\R^n\times\R^m\to\R^n$ by
\begin{equation}\label{Psi}
\Psi(x,v):=f(x)+\nabla\Phi(x)^*v,\quad(x,v)\in\R^n\times\R^m.
\end{equation}
Consider a point $\ox\in\R^n$ satisfying the {\em stationarity condition}
\begin{equation}\label{stat}
0\in f(\ox)+\sub(\th\circ\Phi)(\ox)
\end{equation}
and define the set of {\em Lagrange multipliers} associated with $\ox$ by
\begin{equation}\label{laset}
\Lambda(\ox):=\big\{v\in\R^m\big|\;\Psi(\ox,v)=0,\;v\in\sub\th(\oz)\big\}\;\mbox{ where }\;\oz:=\Phi(\ox).
\end{equation}
We suppose in what follows that $\Lm(\ox)\ne\emp$, which can be ensured under certain qualification conditions discussed in Remark~\ref{socs4}. Observing that it may not hold in the general setting of \eqref{laset} with $\Psi$ from \eqref{Psi}, while sufficient conditions for the existence of Lagrange multipliers in special classes of variational (KKT-type) systems are well known; see Remark~\ref{socs4}. \vspace*{0.03in}

The following basic definition involves the construction $(D\partial\th)(\ox,\oz)\colon\R^{m+1}\tto\R^{m+1}$, i.e., the graphical derivative \eqref{gder} of the first-order subdifferential mapping, which is therefore a {\em second-order} generalized differential construction for $\th\colon\R^m\to\oR$. We then present several equivalent descriptions and calculations for the general case of $\th\in CPWL$.\vspace*{-0.05in}

\begin{Definition}{\bf(critical and noncritical multipliers).}\label{crit} Let $\ox$ satisfy \eqref{stat} with $\Psi$ taken from \eqref{Psi}. Assume that $f$ is differentiable while $\Phi$ is twice differentiable at $\ox$. Then the multiplier $\ov\in\Lambda(\ox)$ is {\sc critical} for \eqref{VS} if there is $0\ne\xi\in\R^n$ satisfying the generalized KKT system
\begin{equation}\label{crc}
0\in\nabla_{x}\Psi(\ox,\ov)\xi+\nabla\Phi(\ox)^*\big(D\sub\th\big)(\oz,\ov)\big(\nabla\Phi(\ox)\xi\big)\;\mbox{ with }\;\oz=\Phi(\ox).
\end{equation}
The multiplier $\ov\in\Lambda(\ox)$ is {\sc noncritical} for \eqref{VS} otherwise, i.e., when the generalized equation \eqref{crc} admits only the trivial solution $\xi=0$.
\end{Definition}\vspace*{-0.1in}

It follows from the calculations below that, in the case where $\th$ is the indicator function of the polyhedral set $\O:=\R^l\times\R^{m-l}_-$ with $0\le l\le m$, our Definition~\ref{crit} reduces to the notions in \cite[Definition~1.41]{is14}, which were introduced by Izmailov \cite{iz05} for pure equality constraints in NLPs and then extended by Izmailov and Solodov \cite{is12} to problems with inequalities. The main advantage of our new setting is that we can efficiently calculate the construction $D\sub\th$ in \eqref{crc} for the general class of CPWL functions, which allows us to deal with a variety of important variational systems appearing in optimization theory and applications.

For any fixed function $\th\in CPWL$ we proceed as follows. Pick a subgradient $\ov\in\sub\th(\oz)$ and introduce the {\em critical cone} for $\th$ at $(\oz,\ov)$ by
\begin{equation}\label{cc}
{\cal K}(\oz,\ov):=\big\{w\in T(\oz;\dom\th)\big|\;\la\bar v,w\ra=d\th(\oz)(w)\big\},
\end{equation}
where the (Dini-Hadamard) {\em subderivative} function $d\th(\oz)\colon\R^n\to\oR$ is defined by
$$
d\th(\oz)(w):=\liminf_{\substack{u\to w\\t\dn 0}}{\frac{\th(\oz+tu)-\th(\oz)}{t}}.
$$
It is shown in \cite[Proposition~10.21]{rw} that for $\th\in CPWL$ the above subderivative construction reduces to the classical directional derivative
\begin{equation}\label{sder}
d\th(\oz)(w)=\th'(\oz;w)=\lim_{t\dn 0}{\frac{\th(\oz+tw)-\th(\oz)}{t}}.
\end{equation}

The critical cone \eqref{cc} agrees with the standard critical cone notion for convex polyhedra; see, e.g., \cite{dr}. Indeed, for $\th=\dd(x;\O)$ we have $\dom\th=\O$ and $d\th(\oz)(w)=0$ for any $w\in T(\oz;\O)$. Thus ${\cal K}(\oz,\ov)=T(\oz;\O)\cap\ov^\bot$ in \eqref{cc}. To avoid confusion, note that in the case of constraint systems in nonlinear programming described by $\Gamma:=\{x\in\R^n|\;g(x)\in\Th\}$ with smooth mappings $g\colon\R^n\to\R^m$, the conventional critical cone as in \cite{bs,is14} is given not in terms of the tangent cone $T(\oz;\Th)$ but via its {\em linearized cone} version
\begin{equation*}
T^{\rm lin}(\ox;\Gamma):=\big\{w\in\R^m\big|\;\nabla g(\ox)w\in T(\oz;\Th)\big\}\;\mbox{ with }\;\oz:=g(\ox).
\end{equation*}\vspace*{-0.15in}

The next proposition calculates the critical cone \eqref{cc} for any function $\th\in CPWL$ in terms of its given data from \eqref{dom} and \eqref{cpwl1}.\vspace*{-0.1in}

\begin{Proposition}{\bf(calculation of the critical cone for CPWL functions).}\label{ccid} Let $\th\in CPWL$ with $(\oz,\ov)\in\gph\partial\th$, and let $\ov_1$, $\ov_2$ from \eqref{eq06} be such that $\ov=\ov_1+\ov_2$. Denote by $K:=K(\oz)$, $I:=I(\oz)$, $J_{1}:=J_+(\oz,\ov_1)$, and $J_{2}:=J_+(\oz,\ov_2)$ the index sets from \eqref{active2} and \eqref{eq05}, respectively. Then the critical cone ${\cal K}(\oz,\ov)$ in \eqref{cc} is calculated by
\begin{equation}\label{cc2}
\begin{array}{ll}
{\cal K}(\oz,\ov)=\Big\{u\in\R^m\Big|&\la a_i-a_j,u\ra=0\;\mbox{ if }\;i,j\in J_1,\\
&\la a_i-a_j,u\ra\le 0\;\mbox{ if }\;(i,j)\in(K\setminus J_1)\times J_1,\\
&\la d_i,u\ra=0\;\mbox{ if }\;i\in J_2,\;\mbox{ and }\;\la d_i,u\ra\le 0\;\mbox{ if }\;i\in I\setminus J_2\;\Big\},
\end{array}
\end{equation}
which means that ${\cal K}(\oz,\ov)={\cal G}_{\tiny\{K,J_1\},\{I,J_2\}}$, where the latter set is defined in \eqref{eq081}.
\end{Proposition}\vspace*{-0.05in}
{\bf Proof}. Picking $u\in{\cal K}(\oz,\ov)$, we show first that $\la a_j,u\ra=\la a_i,u\ra$ whenever $i,j\in J_1$. Taking into account that $\dom d\th(\oz)=T(\oz;\dom\th)$ by \cite[Theorem~10.21]{rw} gives us sequences $t_k\to 0$ and $u_k\to u$ such that $\oz+t_k u_k\in\dom\th$. Thus by passing to a subsequence if necessary, we get a constant index subset $P\subset K$ with $K(\oz+t_k u_k)=P$ for all $k$. It follows from (\ref{sder}) that
\begin{equation}\label{re11}
d\th(\oz)(u)=\la a_s,u\ra\;\mbox{ whenever }\;s\in P.
\end{equation}
If $i\in K$ and $s\in P$, then (\ref{pwlr1}) tells us that $\la a_i,\oz+t_k u_k\ra-\al_i\le\la a_s,\oz+t_k u_k\ra-\al_s$, and so
\begin{equation}\label{re12}
\la a_i,u\ra\le\la a_s,u\ra\;\mbox{ for all }\;i\in K,\;s\in P.
\end{equation}
Furthermore, it follows from (\ref{tanc}) and the choice of $u\in T(\oz;\dom\th)$ that $\la\ov_2,u\ra\le 0$. Employing this together with (\ref{re11}) and (\ref{re12}) gives us the relationships
\begin{equation}\label{re13}
\la a_s,u\ra=d\th(\oz)(u)=\la\ov,u\ra\le\la\ov_1,u\ra=\sum_{i\in J_1}\bar\lm_i\la a_i,u\ra\le\sum_{i\in J_1}\bar\lm_i\la a_s,u\ra=\la a_s,u\ra,
\end{equation}
which yield $\la a_s,u\ra=\sum_{i\in J_1}\bar\lm_i\la a_i,u\ra$. Combining the latter with $\lm_i>0$ for $i\in J_1$ and (\ref{re12}) shows that $\la a_s,u\ra=\la a_i,u\ra$ for any $i\in J_1$ and so $\la a_i-a_j,u\ra=0$ if $i,j\in J_1$.

Consider next the case where $(i,j)\in(K\setminus J_1)\times J_1$. Take $s\in P$ and get from (\ref{re12}) that $\la a_i,u\ra\le\la a_s,u\ra$. Since $\la a_s,u\ra=\la a_j,u\ra$, it tells us that $\la a_i-a_j,u\ra\le 0$. Finally, it follows from (\ref{re13}) that $\la\ov_2,u\ra=0$. Combining this with the inequality
$\la\ov_2,u\ra\le 0$, we arrive at $\la d_i,u\ra=0$ for $i\in J_2$ and $\la d_i,u\ra\le 0$ for $i\in I\setminus J_2$ and thus justify the inclusion ``$\subset$" in (\ref{cc2}).

To prove the opposite inclusion, pick any $u$ from the right-hand side of (\ref{cc2}). It follows from (\ref{tanc}) that $u\in T(\oz;\dom\th)$, which clearly implies that $\la\ov,u\ra\le d\th(\oz)(u)$. Taking into account that $u\in T(\oz;\dom\th)=\dom d\th(\oz)$ and so $d\th(\oz)(u)<\infty$, for any sequence of $\ox+t_tu_k\in\dom d\th(\oz)$ with $t_k\to 0$ and $u_k\to u$ find $P\subset K$ such that $K(\oz+t_k u_k)=P$ for all $k$. Pick further $r\in P$ and observe that  $d\th(\oz)(u)=\la a_r,u\ra$. Then for any $i\in J_1$ we get
$$
\la a_r,u\ra\le\la a_i,u\ra=\la\ov_1,u\ra=\la\ov,u\ra.
$$
This shows that $\la\ov,u\ra=d\th(\oz)(u)$, and hence we arrive at $u\in{\cal K}(\oz,\ov)$, which justifies \eqref{cc2}. $\h$\vspace*{0.05in}

The following theorem provides an equivalent description of critical multipliers from Definition~\ref{crit} for the variational system \eqref{VS} with $\th\in CPWL$ via the critical cone \eqref{cc} calculated in Proposition~\ref{ccid} in terms of the given parameters of $\th$.\vspace*{-0.1in}

\begin{Theorem}{\bf(equivalent description of critical multipliers).}\label{cric002} Let $(\ox,\ov)$ be as in Definition~{\rm\ref{crit}} with $\th\in CPWL$. Then $\ov$ is critical for \eqref{VS} if and only if the primal-dual system
\begin{equation}\label{crc6}
\nabla_{x}\Psi(\ox,\ov)\xi+\nabla\Phi(\ox)^*\eta=0,\;\la\eta,\nabla\Phi(\ox)\xi\ra=0,\;\nabla\Phi(\ox)\xi\in{\cal K}(\oz,\ov),\;\eta\in{\cal K}(\oz,\ov)^*
\end{equation}
with ${\cal K}(\oz,\ov)$ from \eqref{cc2} admits a solution pair $(\xi,\eta)\in\R^n\times\R^m$ with $\xi\ne 0$.
\end{Theorem}\vspace*{-0.05in}
{\bf Proof.} Although we have an independent direct proof of the claimed result, the presented device is based, for brevity, on general facts of Rockafellar's second-order epi-differentiability theory for fully amenable functions; see \cite{roc} and \cite[Chapter~13]{rw}. It follows from \cite[Proposition~13.9]{rw} that $d^2\th(\oz|\ov)(u)= \dd_{{\cal K}(\oz,\ov)}(u)$, where $d^2\th(\oz|\ov)$ stands for the second subderivative of $\th$; see \cite[Definition~13.3]{rw}. Thus we get from  \cite[Theorem~13.40]{rw} that
\begin{equation*}
\big(D\sub\th\big)(\oz,\ov)(u)=\sub\Big(\frac{1}{2}d^2\th(\oz|\ov)\Big)(u)=\sub\dd_{{\cal K}(\oz,\ov)}(u)=N\big(u;{\cal K}(\oz,\ov)\big)
\end{equation*}
for all $u\in\R^m$. Furthermore, it follows from standard convex analysis, Proposition~\ref{ccid}, and the duality relationship in \eqref{eq082} that
\begin{equation*}
N\big(u;{\cal K}(\oz,\ov)\big)={\cal K}(\oz,\ov)^*\cap\{u\}^\bot={\cal F}_{\tiny\{K,J_1\},\{I,J_2\}}\cap\{u\}^\bot,\quad u\in\R^m,
\end{equation*}
where ${\cal F}$ is taken from \eqref{eq080} with the index sets $K,J_1,I,J_2$ defined in Proposition~\ref{ccid}. Comparing this with \eqref{cc2} and Definition~\ref{crit} justifies the claimed statement.$\h$\vspace*{0.05in}

Now we are ready to specify Definition~\ref{crit} in some particular variational systems corresponding to a certain choice of $\th\in CPWL$ therein. Let us start with the original setting of \cite{iz05} for NLPs with pure equality constraints given by ${\cal C}^2$-smooth functions and then proceed with smooth inequality constraints as in \cite{is12,is14}.\vspace*{-0.1in}

\begin{Example}{\bf(critical multipliers in NLPs with equality constraints).}\label{nlpeq} {\rm Choosing the function $\th=\dd_{\{0\}^m}$ in \eqref{VS}, we see that the critical cone \eqref{cc} in this case is $\{0\}^m$, and thus the conditions in \eqref{crc6} are written in the form
$$
\nabla_x\Psi(\ox,\ov)\xi\in\rge\nabla\Phi(\ox)^*,\;\nabla\Phi(\ox)\xi=0
$$
via the range of the adjoint Jacobian. It gives us the definition of critical multipliers in \cite{iz05}.}
\end{Example}\vspace*{-0.1in}

\begin{Example}{\bf(critical multipliers in NLPs with inequality constraints)}.\label{nlpineq} {\rm This case corresponds to $\th=\dd_{\R_-^m}$ in \eqref{VS}. Denote $\Phi=(\ph_1,\ldots,\ph_m)$ and suppose without loss of generality that $I(\oz)=\{1,\ldots,m\}$. For $\ov=(\ov_1,\ldots,\ov_m)\in\sub\th(\oz)=N_{\R_-^m}(\oz)=\R_+^m$ consider the index subsets
$$
I_+(\ov):=\big\{i\in\{1,\ldots,m\}\big|\;\ov_i>0\big\},\;I_0(\ov):=\big\{i\in\{1,\ldots,m\}\big|\;\ov_i=0\big\}
$$
and readily get the critical cone representation
\begin{equation*}
{\cal K}(\oz,\ov)=\big\{u=(u_1,\ldots,u_m)\in\R^m\big|\;u_i=0\;\mbox{ if }\;i\in I_+(\ov)\;\mbox{ and }\;u_i\le 0\;\mbox{ if }\;i\in I_0(\ov)\big\}.
\end{equation*}
Hence conditions \eqref{crc6} read in this case as follows:
\begin{equation}\label{forex1}
\left\{\begin{array}{ll}
\nabla_{x}\Psi(\ox,\ov)\xi+\nabla\Phi(\ox)^*\eta=0,\;\eta=(\eta_1,\ldots,\eta_m)\in\R_+^m,\;\eta_i\nabla\ph_i(\ox)\xi=0\;\mbox{ if }\;i\in I_0(\ov),\\
\nabla\ph_i(\ox)\xi=0\;\mbox{ if }\;i\in I_+(\ov),\;\mbox{ and }\;\nabla\ph_i(\ox)\xi\le 0\;\mbox{ if }\;i\in I_0(\ov),
\end{array}\right.
\end{equation}
which therefore give us the notion of criticality from \cite[Definition~1.41]{is14}.}
\end{Example}\vspace*{-0.1in}

The general case of {\em smooth} equality and inequality constraint systems studied in \cite{is14} is a direct combination of Examples~\ref{nlpeq} and \ref{nlpineq}. In contrast, the following example concerns {\em nonsmooth} constraint systems, where $\th$ is taken as the pointwise maximum function. Such descriptions are particularly appeared in modeling constrained {\em minimax} problems of optimization (see, e.g., \cite{ms152}) and are not covered by the framework of \cite{is14}.\vspace*{-0.05in}

\begin{Example}{\bf(critical multipliers in nonsmooth constraint systems).}\label{minimax} {\rm Consider the variational system \eqref{VS} with $\th(z):=\max\{z_1,\ldots,z_m\}$ for $z=(z_1,\ldots,z_m)\in\R^m$. This function $\th$ is clearly CPWL while nondifferentiable. Taking $(\ox,\ov)$ as in Theorem~\ref{cric002} with $\ov=(\ov_1,\ldots,\ov_m)\in\sub\th(\oz)$ and $\Phi=(\ph_1,\ldots,\ph_m)$, we readily have
\begin{equation*}
\sum_{i\in K(\oz)}\ov_i=1,\;\ov_i\ge 0\;\mbox{ for }\;i\in K(\oz)\;\mbox{ and }\;\ov_i=0\;\mbox{ for }\;i\in \{1,\ldots,m\}\setminus K(\oz),
\end{equation*}
where the index set $K(\oz)$ is taken from (\ref{active2}) and admits the simplification
\begin{equation*}
K:=K(\oz)=\big\{i\in\{1,\ldots,m\}\big|\;\th(\oz)=\oz_i\;\big\}\;\mbox{ with }\;\oz=(\oz_1,\ldots,\oz_m).
\end{equation*}
The critical cone ${\cal K}(\oz,\ov)$ in this framework is represented by
\begin{equation*}
{\cal K}(\oz,\ov)=\big\{u=(u_1,\ldots,u_m)\in\R^m\big|\;\exists\,c\in\R\;\mbox{with}\;u_i=c\;\mbox{if}\;i\in J_1\;\mbox{ and }\;u_i\le c\;\mbox{if}\;i\in K\setminus J_1\big\},
\end{equation*}
where $J_1:=J_+(\oz,\ov)$ is defined in (\ref{eq05}). The direct calculation gives us the dual cone expression
$$
\begin{array}{ll}
{\cal K}(\oz,\ov)^*=\Big\{w=(w_1,\ldots,w_m)\in\R^m\Big|&\disp\sum_{i\in K}w_i=0,\;w_i\ge 0\;\mbox{ if }\;i\in K\setminus J_1\;\mbox{ and }\\
&w_i=0\;\mbox{ if }\;i\in\{1,\ldots,m\}\setminus K\Big\}.
\end{array}
$$
Using the above representations together with Theorem~\ref{cric002} tells us that the criticality of $\ov\in\Lambda(\ox)$ is equivalent to the existence of a solution pair $(\xi,\eta)\in\R^n\times\R^m$, $\xi\ne 0$, of the system
$$
\begin{array}{ll}
\nabla_{x}\Psi(\ox,\ov)\xi+\nabla\Phi(\ox)^*\eta=0,\;\exists\,c\in\R\;\;\mbox{with}\;\nabla\ph_i(\ox)\xi=c\;\mbox{if}\;i\in J_1\;\mbox{ and }\;\nabla\ph_i(\ox)\xi\le c\;\mbox{if}\;i\in K\setminus J_1,\\
\disp\sum_{i\in K}\eta_i=0,\;\eta_i\ge 0\;\mbox{if}\;i\in K\setminus J_1,\;\eta_i=0\;\mbox{if}\;i\in\{1,\ldots,m\}\setminus K,\;\eta_i\big(\nabla\ph_i(\ox)\xi-c\big)=0\;\mbox{if}\;i\in K\setminus J_1,
\end{array}
$$
which provides an explicit construction of critical multipliers in the nonsmooth constraint setting.}
\end{Example}\vspace*{-0.05in}

Finally in this section, we establish an equivalent {\em coderivative} description of critical multipliers in \eqref{VS} with $\th\in CPWL$, which has the potential to be extended beyond the CPWL class and also to problems in infinite-dimensional spaces. Recall first the following relationship between the graphical derivative \eqref{gder} and limiting coderivative \eqref{2.8} of the subdifferential mapping $\partial\th$ established in \cite{rz} and \cite[Theorem~13.57]{rw} for a rather general class of continuously prox-regular and twice epi-differentiable functions including $\th\in CPWL$:
\begin{equation}\label{der-cod}
\big(D\sub\th\big)(\oz,\ov)(u)\subset\big(D^*\sub\th\big)(\oz,\ov)(u),\quad u\in\R^m.
\end{equation}
Furthermore, it is revealed in \cite{rz} that the inclusion in \eqref{der-cod} may be {\em strict} even for smooth functions $\th$ with Lipschitz continuous derivatives.

The next result of its certain independent interest shows that a counterpart of \eqref{der-cod}, with replacing $D^*$ by the {\em regular coderivative} \eqref{2.7} and selecting an appropriate subset of it, holds as {\em equality} at least in the case of $\th\in CPWL$; in fact, in more general settings; see the proof below. This leads us, in particular, to the aforementioned description of critical multipliers.\vspace*{-0.05in}

\begin{Theorem}{\bf(graphical derivative and regular coderivative of the subdifferential mapping for CPWL functions).}\label{gdp10}
Let $\th\in CPWL$ with $(\oz,\ov)\in\gph\sub\th$. Then
\begin{equation}\label{der-cod1}
\dom D\big(\sub\th\big)(\oz,\ov)=-\dom\big(\Hat D^*\sub\th\big)(\oz,\ov)={\cal K}(\oz,\ov)
\end{equation}
via the critical cone ${\cal K}(\oz,\ov)$ calculated in \eqref{cc2}. Moreover, for any $u\in{\cal K}(\oz,\ov)$ we have
\begin{equation}\label{gpco}
\begin{array}{lll}
\big(D\sub\th\big)(\oz,\ov)(u)&=&\big\{w\in\big(\Hat D^*\sub\th\big)(\oz,\ov)(-u)\big|\;\la w,u\ra=0\big\}\\
&=&\Big\{w\Big|\;w\in{\rm{argmin}}\big\{-\la w,u\ra\big|\;w\in\big(\Hat D^*\sub\th\big)(\oz,\ov)(-u)\big\}\Big\}.
\end{array}
\end{equation}
\end{Theorem}\vspace*{-0.05in}
{\bf Proof.} It follows from \cite[Theorem~4.3]{ms15} that
\begin{equation}\label{calpc}
\big(\Hat D^*\sub\th\big)(\oz,\ov)(u)={\cal K}(\oz,\ov)^*\;\mbox{ for any }\;u\in\dom\big(\Hat D^*\sub\th\big)(\oz,\ov)=-{\cal K}(\oz,\ov),
\end{equation}
Also we have from the proof of Theorem~\ref{cric002} above that
\begin{equation*}
\big(D\sub\th\big)(\oz,\ov)(u)=N\big(u;{\cal K}(\oz,\ov)\big)\;\mbox{ whenever }\;u\in{\cal K}(\oz,\ov).
\end{equation*}
Combining these representations shows that \eqref{der-cod1} and the first equality in \eqref{gpco} are satisfied. To verify the second equality in (\ref{gpco}), pick $u\in{\cal K}(\oz,\ov)$ and $w\in(\Hat D^*\sub\th)(\oz,\ov)(-u)$ and then deduce from the maximal monotonicity of $\partial\th$ and \cite[Theorem~2.1]{pr98} that $-\la w,u\ra\ge 0$. Since we always have $0\in(\Hat D^*\sub\th)(\oz,\ov)(-u)$ for any $u\in{\cal K}(\oz,\ov)$, this tells us that
$$
\min\big\{-\la w,u\ra\big|\;w\in\big(\Hat D^*\sub\th\big)(\oz,\ov)(-u)\big\}=0,
$$
which in turn implies the equality
$$
{\rm{argmin}}\big\{-\la w,u\ra\big|\;w\in\big(\Hat D^*\sub\th\big)(\oz,\ov)(-u)\big\}=\big\{w\in\big(\Hat D^*\sub\th\big)(\oz,\ov)(-u)\big|\;\la w,u\ra=0\big\}
$$
and thus completes the proof of the theorem. $\h$\vspace*{0.03in}

It is worth mentioning that the relationships in (\ref{gpco}) between the graphical derivative and regular coderivative of $\sub\th$ can be extended to a more general class of fully amenable functions $\th$ in the sense of \cite[Definition~10.23]{rw}; however, this is beyond the scope of this paper. Let us show now that Theorem~\ref{gdp10} implies the following description of critical multipliers for $\th\in CPWL$.\vspace*{-0.05in}

\begin{Corollary}{\bf(coderivative description of critical multipliers).}\label{cric02} Let $(\ox,\ov)$ be in the setting of Definition~{\rm\ref{crit}} with $\th\in CPWL$. Then $\ov\in\Lambda(\ox)$ is critical for \eqref{VS} if and only if there exists a pair $(\xi,\eta)\in \R^n\times \R^m$ with $\xi\ne 0$ for which
\begin{equation}\label{cric02a}
\eta\in\big(\Hat D^*\sub\th\big)(\ox,\ov)\big(-\nabla\Phi(\ox)\xi\big),\;\nabla_x\Psi(\ox,\ov)\xi+\nabla\Phi(\ox)^*\eta=0,\;\la\eta,\nabla\Phi(\ox)\xi\ra=0.
\end{equation}
\end{Corollary}\vspace*{-0.1in}
{\bf Proof.} The characterizations of critical multipliers in \eqref{cric02a} follow directly from Definition~\ref{crit} and Theorems~\ref{gdp10}, \ref{cric002}. Note also that the construction $\Hat D^*\sub\th$ is calculated in \cite[Theorem~4.3]{ms15} via the given data of $\th\in CPWL$, and these calculations give us efficient descriptions of critical multipliers equivalent to those in \eqref{crc6} with ${\cal K}(\oz,\ov)$ calculated in \eqref{cc2}. $\h$\vspace*{-0.15in}

\section{Noncritical Multipliers and Canonical Perturbations}\sce\vspace*{-0.05in}

This section is devoted to characterizing {\em noncritical} multipliers for \eqref{VS} via some error bound estimating the distance to the solution map of the generalized KKT system as well as via a certain calmness property of the solution map to the {\em canonical perturbation} of \eqref{VS} given by:
\begin{equation}\label{geq2}
\left[\begin{array}{c}
p_1\\p_2
\end{array}
\right]
\in\left[\begin{array}{c}
\Psi(x,v)\\
-\Phi(x)
\end{array}
\right]+\left[\begin{array}{c}
0\\
(\sub\th)^{-1}(v)
\end{array}
\right]
\end{equation}
with the canonical parameter pair $(p_1,p_2)\in\R^n\times\R^m$ and $\th\in CPWL$, where $\Psi$ is defined in \eqref{Psi}. Note that the calmness property, labeled here as {\em semi-isolated calmness}, is different from the conventional calmness and isolated calmness notions for set-valued mappings; see the discussion in Remark~\ref{rem07}.  The next result is an extension of \cite[Proposition~1.43]{is14}, which addresses the case where $\th(\cdot)=\dd(\cdot;\O)$ is the indicator function of the polyhedral set  $\O:=\R^s\times\R^{m-s}_-$, $0\le s\le m-s$, i.e., the classical case of nonlinear programs with $s$ equality and $m-s$ inequality constraints.

Consider the set-valued mapping $G\colon\R^n\times\R^m\tto\R^n\times\R^m$ associated with \eqref{geq2} by
\begin{equation}\label{F-KKT}
G(x,v):=\left[\begin{array}{c}
\Psi(x,v)\\
-\Phi(x)
\end{array}
\right]+\left[\begin{array}{c}
0\\
(\sub\th)^{-1}(v)
\end{array}
\right]
\end{equation}
and then define the {\em solution map} $S\colon\R^n\times\R^n\tto\R^n\times\R^m$ to \eqref{geq2} as the inverse to \eqref{F-KKT} by
\begin{equation}\label{mapS}
S(p_1,p_2):=\big\{(x,v)\in\R^n\times\R^m\big|\;(p_1,p_2)\in G(x,v)\big\}.
\end{equation}\vspace*{-0.35in}

\begin{Theorem}{\bf(characterization of noncritical multipliers via error bound and semi-isolated calmness of solutions under canonical perturbations ).}\label{uplip} Let $(\ox,\ov)\in S(\op_1,\op_2)$ with $(\bar p_1,\bar p_2)=(0,0)$ in \eqref{mapS} under the assumptions of Definition~{\rm\ref{crit}}, where $\th\in CPWL$. Then the following properties of \eqref{VS} and its perturbation \eqref{geq2} are equivalent:

{\bf(i)} The Lagrange multiplier $\ov\in\Lambda(\ox)$ from \eqref{laset} is noncritical for \eqref{VS}.

{\bf(ii)} {\sc (semi-isolated calmness)} There are numbers $\ve>0$, $\ell\ge 0$ and neighborhoods $U$ of $0\in\R^n$ and $W$ of $0\in\R^m$ such that for any $(p_1,p_2)\in U\times W$ and any $(x_{p_1p_2},v_{p_1p_2})\in S(p_1,p_2)\cap\B_\ve(\ox,\ov)$ we have the estimate
\begin{equation}\label{upper}
\|x_{p_1p_2}-\ox\|+\dist\big(v_{p_1p_2};\Lm(\ox)\big)\le\ell\big(\|p_1\|+\|p_2\|\big).
\end{equation}

{\bf(iii)} {\sc (error bound)} There are numbers $\ve>0$  and $\ell\ge 0$ such that the estimate
\begin{equation}\label{subr}
\|x-\ox\|+\dist\big(v;\Lm(\ox)\big)\le\ell\big(\|\Psi(x,v)\|+\dist\big(\Phi(x),\sub\theta^{*}(v)\big)\big)
\end{equation}
holds for any $(x,v)\in \B_\ve(\ox,\ov)$, where  $^*$ signifies the Fenchel conjugate of convex analysis.
\end{Theorem}
{\bf Proof}. To justify first the implication (ii)$\Longrightarrow$(i), we need to verify by Theorem~\ref{cric002} that the validity of  \eqref{upper} ensures that for any solution pair $(\xi,\eta)\in\R^n\times\R^m$ to (\ref{crc6}) we have $\xi=0$. Pick any pair $(\xi,\eta)\in\R^n\times\R^m$  satisfying (\ref{crc6}), let $t>0$, and define $(x_t,v_t):=(\ox+t\xi,\ov+t\eta)$. Thus we have for all $t$ sufficiently small that
$$
\begin{array}{lll}
\Psi(x_t,v_t)-\Psi(\ox,\ov)&=&\big(f(x_t)-f(\ox)\big)+\big(\nabla\Phi(x_t)-\nabla\Phi(\ox)\big)^*\ov+t\nabla\Phi(x_t)^*\eta\\
&=&t\nabla f(\ox)\xi+o(t)+t\big(\nabla^2\Phi(\ox)\xi\big)^*\ov+t\nabla\Phi(\ox)^*\eta+o(t)\\
&=&t\big(\nabla_{x}\Psi(\ox,\ov)\xi+\nabla\Phi(\ox)^*\eta\big)+o(t)=o(t)\;\mbox{ as }\;t\dn 0.
\end{array}
$$
Since $\Psi(\ox,\ov)=0$, we deduce from the last equality that
\begin{equation}\label{eq2}
\Psi(x_t,v_t)=p_{1t}\;\mbox{ with }\;p_{1t}=o(t)\;\mbox{ as }\;t\dn 0.
\end{equation}
Remembering that $\Phi(x_t)=\Phi(\ox)+t\nabla\Phi(\ox)\xi+o(t)$ and letting $z_t:=\Phi(\ox)+t\nabla\Phi(\ox)\xi$ yield
\begin{equation}\label{eq3}
z_t=\Phi(x_t)+p_{2t}\;\mbox{ with }\;p_{2t}=o(t)\;\mbox{ as }\;t\dn 0.
\end{equation}
It is easy to see that $z_t\in\dom\th$ for small $t$, where the set $\dom\th$ is taken from (\ref{dom}).

In what follows we use the notation of Section~2 while denoting for simplicity by $K:=K(\oz)$, $I:=I(\oz)$, $J_{1}:=J_+(\oz,\ov_1)$, and $J_{2}:=J_+(\oz,\ov_2)$ the index sets from {\rm(\ref{active2})} and {\rm(\ref{eq05})}, respectively, with $\oz=\Phi(\ox)$. We proceed with verifying the following statement.\vspace{0.05in}\\
{\bf Claim:} {\em Given $z_t$ as defined above, we have $J_1\subset K(z_t)$ and $J_2\subset I(z_t)$ for all small $t>0$}.\\[1ex]
Starting with checking the inclusion $J_2\subset I(z_t)$, take $i\in J_2$ and get from the definitions that
$$
\la d_i,z_t\ra=\la d_i,\oz\ra+t\la d_i,\nabla\Phi(\ox)\xi\ra=0
$$
due to $\nabla\Phi(\ox)\xi\in{\cal K}(\oz,\ov)$ and $J_2\subset I(\oz)$; thus the second inclusion in the claim holds. To verify the first inclusion therein, pick $i\in J_1$ and check that $z_t\in C_i$, where the polyhedral set $C_i$ is taken from (\ref{pwlr1}). To see this, take $r\in K$ and then get $\la a_r-a_i,\oz\ra=\al_r-\al_i$. It follows from $\nabla\Phi(\ox)\xi\in{\cal K}(\oz,\ov)$ that $\la a_r-a_i,\nabla\Phi(\ox)\xi\ra\le 0$. These lead us to
$\la a_r-a_i,z_t\ra\le\al_r-\al_i$ for $r\in K$. Similarly we can show that $\la a_r-a_i,z_t\ra\le\al_r-\al_i$ for $r\in T_1\setminus K$. Thus we arrive at $\la a_r-a_i,z_t\ra\le\al_r-\al_i$ for $r\in T_1$, and hence $z_t\in C_i$ while completing the proof of this Claim.\vspace*{0.05in}

Let us next show that $v_t\in\sub\th(z_t)$ whenever $t$ is sufficiently small. Indeed, it follows from $\eta\in({{\cal K}(\oz,\ov)})^*$, \eqref{eq082}, and Proposition~\ref{ccid} that
$$
\eta\in{\cal K}(\oz,\ov)^*={\cal G}^{*}_{\tiny\{K,J_1\},\{I,J_2\}}={\cal F}_{\tiny\{K,J_1\},\{I,J_2\}}.
$$
Then using \eqref{eq080} gives us the representation $\eta=\eta_1+\eta_2$ such that
\begin{equation*}
\begin{array}{lll}
\displaystyle{\eta_1:=\sum_{i,j\in J_1}{\beta_{ij}(a_i-a_j)}+\sum_{(i,j)\in(K\setminus J_1)\times J_1}{\rho_{ij}(a_i-a_j)}}&\mbox{and}&\disp{\eta_2:={\sum_{s\in J_2}\tau_{1s}d_s+\sum_{s\in I\setminus J_2}\tau_{2s}d_s},}\\
\beta_{ij}\in\R\quad\mbox{for}\quad i,j\in J_1&\mbox{and}&\tau_{1s}\in\R\quad\mbox{for}\quad s\in J_2,\\
\rho_{ij}\ge 0\quad\mbox{for}\quad(i,j)\in(K\setminus J_1)\times J_1 &\mbox{and}&\tau_{2s}\ge 0\quad\mbox{for}\quad s\in I\setminus J_2.
\end{array}
\end{equation*}
We know that $K(z_t)\subset K(\oz)$ and $I(z_t)\subset I(\oz)$ whenever $t$ is small enough. Picking $i_0\in K(\oz)\setminus K(z_t)$ and $j\in J_1$, deduce from the above Claim that $j\in K(z_t)$, which together with $\nabla\Phi(\ox)\xi\in{\cal K}(\oz,\ov)$ brings us to $\la a_{i_0}-a_j,\nabla \Phi(\ox)\xi\ra<0$.
This implies by $\la\eta,\nabla\Phi(\ox)\xi\ra=0$ that $\rho_{i_0j}=0$ in the expression of $\eta_1$. Thus we arrive by using (\ref{eq06}) at the relationships
\begin{equation}\label{pnc22}
\begin{array}{lll}
v_{1t}:=\ov_1+t\eta_1&=&\disp{\sum_{i\in K}\bar\lm_i a_i+t\sum_{i,j\in J_1}{\beta_{ij}(a_i-a_j)}+t\sum_{(i,j)\in(K\setminus J_1)\times J_1}{\rho_{ij}(a_i -a_j)}}\\
&=&\disp{\sum_{i\in J_1}\bar\lm_i a_i+ t\sum_{i,j\in J_1}{\beta_{ij}(a_i-a_j)}+t\sum_{(i,j)\in(K(z_k)\setminus J_1)\times J_1}{\rho_{ij}(a_i-a_j)}}.
\end{array}
\end{equation}
When $t$ is small, there are $\lm'_{ti}\ge 0$ for $i\in K(z_t)$ such that $\sum_{i\in K(z_t)}\lm'_{ti}=\sum_{i\in K}\bar\lm_i=1$ and
\begin{equation}\label{pnc2}
v_{1t}=\sum_{i\in J_1}\lm'_{ti} a_i+\sum_{i\in K(z_t)\setminus J_1 }\lm'_{ti}a_i.
\end{equation}
Similarly, pick $s_0\in I(\oz)\setminus I(z_t)$ and observe by $\nabla\Phi(\ox)\xi\in{\cal K}(\oz,\ov)$ that $\la d_{s_0},\nabla\Phi(\ox)\xi\ra<0$. Thus we get from $\la\eta,\nabla\Phi(\ox)\xi\ra=0$ that $\tau_{2s_0}=0$ above, which ensures in turn that
\begin{equation}\label{pnc3}
\begin{array}{lll}
v_{2t}:=\ov_2+t\eta_2&=&\disp{\sum_{s\in I}\bar\mu_s d_s+t{\sum_{s\in J_2}\tau_{1s} d_s+t\sum_{s\in I\setminus J_2}\tau_{2s}d_s}}\\
&=&\disp{\sum_{s\in I\setminus J_2}(\bar\mu_s+t\tau_{2s})d_s+\sum_{s\in J_2}(\bar\mu_s+t\tau_{1s})d_s}\\
&=&\disp{\sum_{s\in I(z_t)\setminus J_2}(t\tau_{2s})d_s+\sum_{s\in J_2}(\bar\mu_s+t\tau_{1s})d_s}.
\end{array}
\end{equation}
Employing now (\ref{pnc2}) and (\ref{pnc3}) together with the above Claim shows that
$$
v_t=v_{1t}+v_{2t}\in{\rm co}\big\{a_i\big|\;i\in K(z_t)\big\}+\Big\{\sum_{s\in I(z_t)}\mu_s d_s\Big|\;\mu_s\ge 0\Big\}=\sub\th(z_t)
$$
as desired. Using this along with (\ref{eq2}) and (\ref{eq3}) tells us that $(x_t,v_t)$ is a solution to (\ref{geq2}) associated with $(p_{1t},p_{2t})$, and hence we arrive at
$$
t\|\xi\|=\|x_t-\ox\|\le\ell\big(\|p_{1t}\|+\|p_{2t}\|\big)=\ell\|o(t)\|
$$
by \eqref{upper}. It yields $\xi=0$ and thus justifies the claimed implication (ii)$\Longrightarrow$(i).\vspace*{0.05in}

To verify the opposite one (i)$\Longrightarrow$(ii), it suffices to check that under the validity of (i) there are numbers $\ve>0$, $\ell\ge 0$ and neighborhoods $U$ of $0\in\R^n$ and $W$ of $0\in\R^m$ such that for any $(p_1,p_2)\in U\times W$ and any $(x_{p_1p_2},v_{p_1p_2})\in S(p_1,p_2)\cap\B_\ve(\ox,\ov)$ we have the estimate
\begin{equation}\label{upper2}
\|x_{p_1p_2}-\ox\|\le\ell\big(\|p_1\|+\|p_2\|\big).
\end{equation}
Indeed, assuming for the moment that \eqref{upper2} holds and showing then that there is $\ell'\geq 0$ with
\begin{equation}\label{up2}
\dist\big(v_{p_1p_2};\Lm(\ox)\big)\le\ell'\big(\|x_{p_1p_2}-\ox\|+\|p_1\|+\|p_2\|\big),
\end{equation}
we immediately get (\ref{upper}). Let us first justify the validity of \eqref{up2}. To furnish this, observe that the subdifferential $\sub\th(\oz)$ is a convex polyhedral set and find by the classical Minkowski-Weyl theorem $r>0$, $q\in\R^r$, and $A\in\R^{m\times r}$ such that $\sub\th(\oz)$ is represented in the form
$$
\sub\th(\oz)=\big\{y\in\R^m\big|\;Ay\le q\;\big\}.
$$
For any vectors $a\in\R^n$ and $b\in\R^r$, define now the set
\begin{equation}\label{eq4}
{\cal D}_{\ox}(a,b)=\big\{v\in\R^m\big|\;\Psi(\ox,v)=a,\;Av\le b\big\},
\end{equation}
and observe that ${\cal D}_{\ox}(0,q)=\Lm(\ox)$ with $\Lm(\ox)$ given in (\ref{laset}). It follows from \cite[Proposition~3.3(i)]{ms15} that $v_{p_1p_2}\in\sub\th(\Phi(x_{p_1p_2})+p_2)\subset\sub\th(\oz)$ whenever $(p_1,p_2)\in U\times W$. Denoting by $L\ge 0$ a common Lipschitz constant for the mappings $f$, $\nabla\Phi$ and employing the classical Hoffman Lemma, we find a positive constant $M$ such that
\begin{equation}\label{up3}
\begin{array}{lll}
\dist\big(v_{p_1p_2};\Lm(\ox)\big)&=&\dist\big(v_{p_1p_2};{\cal D}_{\ox}(0,q)\big)\le M\|\Psi(\ox,v_{p_1p_2})\|\\\\
&\le&M\big(\|\Psi(\ox,v_{p_1p_2})-\Psi(x_{p_1p_2},v_{p_1p_2})\|+\|\Psi(x_{p_1p_2},v_{p_1p_2})\|\big)\\\\
&\le&M\big(L\|x_{p_1p_2}-\ox\|+\|p_1\|\big)\le M\big(L\|x_{p_1p_2}-\ox\|+\|p_1\|+\|p_2\|\big),
\end{array}
\end{equation}
which justifies (\ref{up2}). To complete the proof of the theorem, it thus remains to verify (\ref{upper2}).

Suppose on the contrary that \eqref{upper2} fails, i.e., for any $k\in\N$ there are $(p_{1k},p_{2k})\in\B_{\frac{1}{k}}(0)\times\B_{\frac{1}{k}}(0)$ and $(x_k,v_k)\in S(p_{1k},p_{2k})\cap\B_{\frac{1}{k}}(\ox,\ov)$ satisfying
$$
\frac{\|x_k-\ox\|}{\|p_{1k}\|+\|p_{2k}\|}\to\infty\Longleftrightarrow\frac{\|p_{1k}\|+\|p_{2k}\|}{\|x_k-\ox\|}\to 0\;\mbox{ as }\;k\to\infty,
$$
which yields $p_{1k}=o(\|x_k-\ox\|)$ and $p_{2k}=o(\|x_k-\ox\|)$. Let $z_k:=\Phi(x_k)+p_{2k}$ and observe by (\ref{geq2}) that $(z_k,v_k)\in\gph\sub\th$. Applying Theorem~\ref{gphpar} tells us that $J_1\subset K(z_k)\subset K(\oz)$ and $J_2\subset I(z_k)\subset I(\oz)$. Passing to a subsequence of $(z_k,v_k)$ if necessary, we suppose without loss of generality that there are subsets $P\subset K(\oz)$ and $Q\subset I(\oz)$ such that
\begin{equation*}
P=K(z_k)\quad\mbox{and}\quad Q=I(z_k)\;\;\mbox{whenever}\;\;k\in\N.
\end{equation*}
Remember that for each $k\in\N$ the pair $(x_k,v_k)$ solves the perturbed system (\ref{geq2}) associated with the parameter pair $(p_{1k},p_{2k})$. Thus we have
\begin{equation}\label{up10}
\begin{array}{lll}
o(\|x_k-\ox\|)=p_{1k}&=&\Psi(x_{k},v_{k})=\Psi(x_{k},\ov)-\Psi(\ox,\ov)+\nabla\Phi(x_k)^*(v_k-\ov)\\
&=&\nabla_{x}\Psi(\ox,\ov)(x_k-\ox)+\nabla\Phi(\ox)^*(v_k-\ov)+o(\|x_k-\ox\|).
\end{array}
\end{equation}
Employing (\ref{fos}) together with $v_k\in\sub\th(z_k)$, we find $\lm_{ik}\ge 0$ with $i\in P$ and $\mu_{ik}\geq 0$ with $i\in Q$ so that each $v_{k}$ is represented in the form $v_k=v_{1k}+v_{2k}$, where
$$
v_{1k}=\sum_{i\in P}\lm_{ik}a_i\;\mbox{ and }\;v_{2k}=\sum_{i\in Q}\mu_{ik}d_i\;\mbox{ with }\;\sum_{i\in P}{\lm_{ik}}=1.
$$
Combining this with (\ref{eq06}) and (\ref{up10}) implies that
\begin{equation}\label{up5}
\begin{array}{lll}
\disp-\nabla_{x}\Psi(\ox,\ov)\frac{(x_k-\ox)}{\|x_k-\ox\|}+\frac{o(\|x_k-\ox\|)}{\|x_k-\ox\|}=\disp\frac{1}{\|x_k-\ox\|}\nabla\Phi(\ox)^*\Big[(v_{1k}-\ov_1)+(v_{2k}-\ov_2)\Big ]\\
~~~~~~~~~=\disp\frac{1}{\|x_k-\ox\|}\nabla\Phi(\ox)^*\Big[\Big(\sum_{i\in P}\lm_{ik}a_i-\sum_{j\in J_1}\bar\lm_ja_j\Big)+\Big(\sum_{i\in Q}\mu_{ik} d_i-\sum_{j\in J_2}\bar\mu_jd_j\Big)\Big]\\
~~~~~~~~~=\disp\frac{1}{\|x_k-\ox\|}\nabla\Phi(\ox)^*\Big[\sum_{i\in P}\lm_{ik}\sum_{j\in J_1}\bar\lm_j(a_i-a_j)+\Big(\sum_{i\in Q}\mu_{ik}d_i-\sum_{j\in J_2}\bar\mu_jd_j\Big)\Big]\\
~~~~~~~~~\in\nabla\Phi(\ox)^*\big(\span\big\{a_i-a_j\big|\;i,j\in J_1\big\}+\cone\big\{a_i-a_j\big|\;(i,j)\in(P\setminus J_1)\times J_1\big\}\\
~~~~~~~~~~~~~~~~~~~~~~+\cone\big\{d_j\big|\;j\in Q\setminus J_2\big\}+\span\big\{d_j\big|\;j\in J_2\big\}\big).
\end{array}
\end{equation}
Assume without loss of generality that
\begin{equation}\label{xi}
\frac{x_k-\ox}{\|x_k-\ox\|}\to\xi\;\mbox{ as }\;k\to\infty\;\mbox{ for some }\;\xi\ne 0.
\end{equation}
Since the set on the right-hand side of (\ref{up5}) is closed, by passing to the limit as $k\to\infty$ we get
$$
\begin{array}{lll}
-\nabla_{x}\Psi(\ox,\ov)\xi\in\nabla\Phi(\ox)^*\big(\span\big\{a_i-a_j\big|\;i,j\in J_1\big\}+\cone\big\{a_i-a_j\big|\;(i,j)\in(P\setminus J_1)\times J_1\big\}\\
~~~~~~~~~~~~~~~~~~~~~~+\cone\big\{d_j\big|\;j\in Q\setminus J_2\big\}+\span\big\{d_j\big|\;j\in J_2\big\}\big ).
\end{array}
$$
This allows us to find a vector $\eta\in\R^m$ in the form $\eta=\eta_1+\eta_2$ with
\begin{equation*}
\begin{array}{lll}
\eta_1\in\span\big\{a_i-a_j\big|\;i,j\in J_1\big\}+\cone\big\{a_i-a_j\big|\;(i,j)\in(P\setminus J_1)\times J_1\big\},\\
\eta_2\in\cone\big\{d_j\;\big|\;j\in Q\setminus J_2\big\}+\span\big\{d_j\big|\;j\in J_2\big\}
\end{array}
\end{equation*}
for which $\nabla_{x}\Psi(\ox,\ov)\xi+\nabla\Phi(\ox)^*\eta=0$, i.e., the first formula in \eqref{crc6} holds. We clearly have
\begin{equation}\label{eq7}
\eta_1=\disp{\sum_{i,j\in J_1}\gg_{ij}(a_i-a_j)+\sum_{(i,j)\in(P\setminus J_1)\times J_1}\gg'_{ij}(a_i-a_j)},\quad\eta_2=\sum_{t\in J_2}\tau_{t}d_t+\sum_{t\in Q\setminus J_2}\tau'_td_t
\end{equation}
with some numbers $\gg_{ij}\in\R$, $\gg'_{ij}\ge 0$, $\tau_t\ge 0$, and $\tau'_t\in\R$. Furthermore, it follows from \eqref{eq080}--\eqref{eq082} by taking into account the inclusions $P\subset K(\oz)$ and $Q\subset I(\oz)$ for the index sets $P$ and $Q$ selected above that $\eta\in{\cal K}(\oz,\ov)^*$, which is the last condition in \eqref{crc6}.

We now claim that $\nabla\Phi(\ox)\xi\in{\cal K}(\oz,\ov)$, which is the third condition in \eqref{crc6}. To verify the claim, pick $i,j\in J_1$ and conclude by the inclusion $J_1\subset P\subset K(\oz)$ together with (\ref{eq3}) that
$$
\la a_i-a_j,z_k-\oz\ra=\la a_i-a_j,\Phi(x_k)+p_{2k}-\Phi(\ox)\ra=0
$$
from where we deduce by $p_{2k}=o(\|x_k-\ox\|)$ the equality
$$
\Big\la a_i-a_j,\nabla\Phi(\ox)\frac{x_k-\ox}{\|x_k-\ox\|}+\frac{o(\|x_k-\ox\|)}{\|x_k-\ox\|}\Big\ra=0.
$$
By passing to the limit therein as $k\to\infty$ with using \eqref{xi}, this results in
\begin{equation}\label{eqq1}
\la a_i-a_j,\nabla\Phi(\ox)\xi\ra=0\;\mbox{ for all }\;i,j\in J_1.
\end{equation}
Taking $i\in K\setminus J_1$, $j\in J_1$ and proceeding similarly to the above, we get
\begin{equation}\label{eqq2}
\la a_i-a_j,\nabla\Phi(\ox)\xi\ra\le 0\;\mbox{ whenever }\;(i,j)\in(K\setminus J_1)\times J_1.
\end{equation}
Pick now $t\in J_2$ and observe that $\la d_t,\Phi(x_k)+p_{2k}-\Phi(\ox)\ra=0$ by the inclusion $J_2\subset Q\subset I$. Combining this with $p_{2k}=o(\|x_k-\ox\|)$ gives us
$$
\Big\la d_t,\nabla\Phi(\ox)\frac{x_k-\ox}{\|x_k-\ox\|}+\frac{o(\|x_k-\ox\|)}{\|x_k-\ox\|}\Big\ra=0,
$$
which allows to establish the equality
\begin{equation}\label{eqq3}
\la d_t,\nabla\Phi(\ox)\xi\ra=0\;\mbox{ whenever }\;t\in J_2.
\end{equation}
Furthermore, for any $t\in I\setminus J_2$ we have $\la d_t,\Phi(x_k)+p_{2k}-\Phi(\ox)\ra\le 0$, which implies by the similar arguments that $\la d_t,\nabla \Phi(\ox)\xi\ra\le 0$. Using this together with (\ref{eqq1})--(\ref{eqq3}) and representation (\ref{cc2}) tells us that $\nabla\Phi(\ox)\xi\in{\cal K}(\oz,\ov)$,
and thus $\nabla\Phi(\ox)\xi\in{\cal K}(\oz,\ov)$.

To employ finally Theorem~\ref{cric002}, it remains to verify the second equality in \eqref{crc6}. It is easy to see that (\ref{eqq1}) holds if $J_1$ is replaced by $P$. Similarly, inequality (\ref{eqq3}) is still true provided that $J_2$ is replaced by $Q$. Using these observations along with (\ref{eq7}), we arrive at $\la\eta,\nabla\Phi(\ox)\xi\ra=0$, which confirms that the pair $(\xi,\eta)$ satisfies all the conditions in (\ref{crc6}). By assertion (i) of the theorem we know that $\ov\in\Lambda(\ox)$ is a noncritical multiplier for \eqref{VS}, and so $\xi=0$ by Theorem~\ref{cric002}, which thus contradicts \eqref{xi}. This justifies (i)$\Longrightarrow$(ii).

Since implication (iii)$\Longrightarrow$(ii) is trivial, it remains to justify implication (ii)$\Longrightarrow$(iii). To this end we first show that there are numbers $\ve>0$, $\dd>0$, and $\ell\ge 0$ such that the estimate
\begin{equation}\label{subr3}
\|x-\ox\|+\dist\big(v;\Lm(\ox)\big)\le\ell\Big(\|\Psi(x,v)\|+\dist\big(0,\big[-\Phi(x)+(\sub\theta)^{-1}(v)\big]\cap\B_\dd\big)\Big)
\end{equation}
holds for any $(x,v)\in \B_\ve(\ox,\ov)$. To proceed, take $\ve$ and the neighborhoods $U$ and $V$ from (ii) and choose $\ve'\le \ve$ so such that for any $(x,v)\in \B_{\ve'}(\ox,\ov)$ we have $\Psi(x,v)\in U $. Pick $(x,v)\in\B_{\ve'}(\ox,\ov)$, let $\dd>0$ with $\B_{\dd}(0)\subset V$, and then  observe that estimate \eqref{subr3} trivially holds if $[-\Phi(x)+(\sub\theta)^{-1}(v)]\cap\B_\dd =\emp$. Thus we can assume that $[-\Phi(x)+(\sub\theta)^{-1}(v)]\cap\B_\dd\ne\emp$. Since the set $[-\Phi(x)+(\sub\theta)^{-1}(v)]\cap\B_\dd$ is closed, there is $p_2\in[-\Phi(x)+(\sub\theta)^{-1}(v)]\cap\B_\dd$ with
$$
\dist\big(0,\big[-\Phi(x)+(\sub\theta)^{-1}(v)\big]\cap\B_\dd \big)=\|p_2\|.
$$
Denote $p_1:=\Psi(x,v)$ and observe that $(p_1,p_2)\in U\times V$, which yields $(x,v)\in S(p_1,p_2)\cap\B_{\ve'}(\ox,\ov)$. Then (ii) gives us the relationships
$$
\|x-\ox\|+\dist\big(v;\Lm(\ox)\big)\le\ell\big(\|p_1\|+\|p_2\|\big)=\ell\Big(\|\Psi(x,v)\|+\dist\big(0,\big[-\Phi(x)+(\sub\theta)^{-1}(v)\big]\cap\B_\dd\big)\Big)
$$
and hence verifies (\ref{subr3}). To derive further estimate (\ref{subr}), pick $\dd'<\dd$ and take $(x,v)\in\B_\dd(\ox,\ov)\cap\B_{\dd'\ell}(\ox,\ov)$. If $[-\Phi(x)+(\sub\theta)^{-1}(v)]\cap\B_{\dd'}\ne\emp$, then 
$$
\dist\big(0,\big[-\Phi(x)+(\sub\theta)^{-1}(v)\big]\cap\B_{\dd'}\big)=\dist\big(0,-\Phi(x)+(\sub\theta)^{-1}(v)\big),
$$
which leads us to the equality
$$
\dist\big(0,\big[-\Phi(x)+(\sub\theta)^{-1}(v)\big]\cap\B_{\dd}\big)=\dist\big(0,-\Phi(x)+(\sub\theta)^{-1}(v)\big).
$$
Otherwise we have $[-\Phi(x)+(\sub\theta)^{-1}(v)]\cap\B_{\dd'}=\emp$ telling us that
$$
\ell\dist\big(0,-\Phi(x)+(\sub\theta)^{-1}(v)\big)\ge\ell\dd'\ge\|x-\ox\|+\|v-\ov\|\ge\|x-\ox\|+\dist\big(v;\Lm(\ox)\big).
$$
Taking finally into account the well-known relationship $(\sub\theta)^{-1}=\sub\th^*$ justifies property (iii) and thus completes the proof of the theorem. $\h$\vspace*{0.05in}

Next we discuss specifications of the error bound \eqref{subr} for particular forms of $\th\in CPWL$ considered in Examples~\ref{nlpeq}, \ref{nlpineq}, and \ref{minimax}.\vspace*{-0.05in}

\begin{Example} {\bf(specifications of error bound).}\label{error-ineq} {\rm We examine the following three cases of the function $\th$ in \eqref{subr} corresponding to the settings of Examples~\ref{nlpeq}, \ref{nlpineq}, and \ref{minimax}, respectively.

{(\bf i)} $\th=\dd_{\{0\}^m}$. This is the case of NLPs with equality constraints, where $\sub\th=N_{\{0\}^m}$ and hence we have the representation
$$
\sub\th^*(v)=N(v;\R^m)=\{0\}\;\mbox{ for all }\;v\in\R^m.
$$
Employing the latter, the error bound (\ref{subr}) reduces to the estimate
\begin{equation}\label{rt1}
\|x-\ox\|+\dist\big(v;\Lm(\ox)\big)\le\ell\big(\|\Psi(x,v)\|+\|\Phi(x)\|\big),
\end{equation}
which is the same as in \cite[Theorem~1.43]{is14} for NLPs with equality constraints.

{(\bf ii)} $\th=\dd_{\R^m_-}$. This reminds us Example~\ref{nlpineq}. In this case we have
$$
\sub\th^*(v)=N_{\R^m_+}(v),
$$
which implies that the error bound (\ref{subr}) is equivalent to the estimate
\begin{equation*}
\|x-\ox\|+\dist\big(v;\Lm(\ox)\big)\le\ell\Big(\|\Psi(x,v)\|+\dist\big(\Phi(x),N(v;\R^m_+)\big)\Big ).
\end{equation*}
Employing now the arguments similar to \cite[Theorem~2]{f02} shows that the latter inequality amounts to the existence of numbers $\dd>0$ and $M\ge 0$ so that for any $(x,v)\in\B_\dd(\ox,\ov)$ the estimate
\begin{equation*}
\|x-\ox\|+\dist\big(v;\Lm(\ox)\big)\le M\big(\|\Psi(x,v)\|+\|\min\{v,-\Phi(x)\}\|\big)
\end{equation*}
holds, which is the well-known error bound property for KKT systems with inequality constraints; see \cite[Theorem~1.43]{is14} and \cite[Proposition~6.2.7]{fp} for more details.

{(\bf iii)} $\th(z)=\max\{z_1,\ldots,z_m\}$ for $z=(z_1,\ldots,z_m)\in\R^m$. It is easy to observe that the function $\th$ can be equivalently written as
$$
\th(z)=\sup_{y\in M}\big\{\la z,y\ra\big\}\;\mbox{ with }\;M:=\Big\{y=(y_1,\ldots,y_m)\Big|\;\sum_{i=1}^{m}y_i=1,\;y_i\ge 0\Big\}.
$$
This readily tells us that
$$
\sub\th^*(v)=N(v;M)\;\mbox{ for all }\;v\in\R^m,
$$
and therefore the error bound property (\ref{subr}) reduces to
$$
\|x-\ox\|+\dist\big(v;\Lm(\ox)\big)\le\ell\Big(\|\Psi(x,v)\|+\dist\big(\Phi(x),N(v;M)\big)\Big).
$$}
\end{Example}

Note that the proof Theorem~\ref{uplip} is heavily based on the {\em second-order subdifferential calculations} for $\th\in CPWL$ conducted in \cite{ms15} being different from the one given \cite[Proposition~1.43]{is14} for the classical KKT system with $\th=\dd(z;\R^s_-\times\R^{m-s})$ in \eqref{geq2}. At the same time, we employ some arguments developed in the proof of the aforementioned result from \cite{is14}.\vspace*{0.03in}

It is not hard to deduce from Theorem~\ref{cric002} (see below) that the condition
\begin{equation}\label{sosc}
\la\nabla_{x}\Psi(\ox,\ov)u,u\ra>0\;\mbox{ for all }\;0\ne u\in\R^n\;\mbox{ with }\;\nabla\Phi(\ox)u\in{\cal K}(\oz,\ov)
\end{equation}
is {\em sufficient} for the multiplier $\ov\in\Lm(\ox)$ to be {\em noncritical}. Consider now its implementation for the problems of {\em composite optimization} formulated in \eqref{comp}, where $\ph_0\colon\R^n\to\R$ and $\Phi\colon\R^n\to\R^m$ are twice differentiable at the reference point, and where $\th\in CPWL$ for $\th\colon\R^m\to\oR$. Note that problem \eqref{comp} can be written in conventional constrained optimization form:
\begin{equation}\label{comp1}
\mbox{minimize }\;\ph_0(x)+\big(\th\circ\Phi\big)(x)\;\mbox{ subject to }\;\Phi(x)\in\dom\th.
\end{equation}
Pick a feasible solution $\ox$ to \eqref{comp} (i.e., such $\ox$ where $\th(\oz)<\infty$ with $\oz=\Phi(\ox)$ and define in terms of the Lagrangian \eqref{lagr} the collection of {\em Lagrange multipliers} for \eqref{comp} at $\ox$ given by
\begin{equation}\label{laset2}
\Lambda_{\small{\rm com}}(\ox):=\big\{v\in\R^m\big|\;\nabla_x L(\ox,v)=0,\;v\in\sub\th(\oz)\big\}\;\mbox{ with }\;\oz:=\Phi(\ox).
\end{equation}
Observe that the set of Lagrange multipliers \eqref{laset} for the general variational system \eqref{VS} studied above reduces to the one in \eqref{laset2} for the composite optimization problem \eqref{comp} by putting $\Psi=\nabla_x L$. In this case the sufficient condition \eqref{sosc} for noncriticality in \eqref{VS} reads as
\begin{equation}\label{sosc2}
\la\nabla_{xx}^2L(\ox,\ov)u,u\ra>0\;\mbox{ for all }\;0\ne u\in\R^n\;\mbox{ with }\;\nabla\Phi(\ox)u\in {\cal K}(\oz,\ov),
\end{equation}
which is a usual form of second-order sufficient conditions for various problems of constrained optimization; see, e.g., \cite{bs,is14}. We show now that \eqref{sosc2} gives us a {\em second-order sufficient condition} (SOSC) for {\em strict minimizers} in the general class \eqref{comp1} under consideration.\vspace*{-0.05in}

\begin{Theorem}{\bf(sufficient condition for strict local minimizers and multiplier noncriticality in composite optimization).}\label{sosc3} Let $\ox$ be a feasible solution to \eqref{comp} such that $\Lambda_{\small{\rm com}}(\ox)\ne\emp$, and let $\ov\in\Lambda_{\small{\rm com}}(\ox)$. Then the validity of \eqref{sosc2} ensures that $\ox$ is a strict local minimizer for \eqref{comp}, i.e., $\ph(\ox)<\ph(x)$ for any $x\ne\ox$ sufficiently close to $\ox$. Furthermore, any $\ov\in\Lambda_{\small{\rm com}}(\ox)$ for which \eqref{sosc2} holds is a noncritical multiplier for \eqref{comp} associated with $\ox$.
\end{Theorem}\vspace*{-0.05in}
{\bf Proof.} Suppose on the contrary that $\ox$ is not a strict local minimizer for \eqref{comp} under the conditions of the theorem. Then we find a sequence of $x_k\to\ox$ as $k\to\infty$ for which we have
\begin{equation*}
\ph_0(x_k)+\th\big(\Phi(x_k)\big)\le\ph_0(\ox)+\th\big(\Phi(\ox)\big)\;\mbox{ and }\;\Phi(x_k)\in\dom\th.
\end{equation*}
Denoting $z_k:=\Phi(x_k)$, we have $K(z_k)\subset K(\oz)$ for all large $k$, where the active index set $K(\cdot)$ is defined in \eqref{active2}.
Extracting a subsequence of $\{z_k\}$ if necessary, find a constant subset $P\subset K(\oz)$ so that $K(z_k)=P$ for all $k$. Define $u_k\:=\frac{x_k-\ox}{\|x_k-\ox\|}$ and suppose without loss of generality that $u_k\to\ou$ as $k\to\infty$ for some $\ou\in\R^n$. Then for any $r\in P$ we have  by the choice of $P$ that
$$
\th\big(\Phi(\ox)\big)=\la a_r,\Phi(\ox)\ra-\al_r\;\mbox{ and }\;\th\big(\Phi(x_k)\big)=\la a_r,\Phi(x_k)\ra-\al_r,
$$
where $(a_r,\al_k)$ are taken from the description of $\th$ in \eqref{eq00}. Hence
\begin{equation}\label{eq31}
\big(\ph_0(x_k)-\ph_0(\ox)\big)+\la a_r,z_k-\oz\ra\le 0,
\end{equation}
which clearly leads us to the inequality
\begin{equation}\label{eq32}
\nabla\ph_0(\ox)\ou+\la a_r,\nabla\Phi(\ox)\ou\ra\le 0\;\mbox{ for any }\;r\in P.
\end{equation}
Invoking $\ov_1,\ov_2$ from \eqref{eq06} and taking into account that $\ov_1\in\sub\th(\oz)$ and $a_r\in\sub\th(z_k)$ whenever $r\in P$, we deduce from the convexity of $\th$ that
$$
\la\ov_1,z_k-\oz\ra\le\la a_r,z_k-\oz\ra,\quad r\in P.
$$
Combining this with (\ref{eq32}) tells us that
\begin{equation}\label{eq33}
\nabla\ph_0(\ox)\ou+\la\ov_1,\nabla\Phi(\ox)\ou\ra\le 0.
\end{equation}
Moreover, by the inclusion $I(z_k)\subset I(\oz)$ we arrive at
\begin{equation}\label{eq34}
\la\ov_2,z_k-\oz\ra\le 0\;\mbox{ and }\;\la\ov_2,\nabla\Phi(\ox)\ou\ra\le 0.
\end{equation}
Since $\ov\in\Lambda_{\small\mbox{com}}(\ox)$, it follows that $\nabla\ph_0(\ox)+\nabla\Phi(\ox)^*\ov=0$, which being combined with
(\ref{eq33}) and (\ref{eq34}) leads us to the equalities
\begin{equation}\label{eq35}
\la a_j,\nabla\Phi(\ox)\ou\ra=-\nabla\ph_0(\ox)\ou\;\mbox{ for }\;j\in J_1\;\mbox{ and }\;\la d_t,\nabla\Phi(\ox)\ou\ra=0\;\mbox{ for }\;t\in J_2.
\end{equation}
Pick now $i\in K:=K(\oz)$ and deduce from the convexity of $\th$ that
$$
\la a_i,z_k-\oz\ra\le\th(z_k)-\th(\oz)\le-\big(\ph_0(x_k)-\ph_0(\ox)\big)
$$
by which we obtain the relationships
$$
\la a_i,\nabla\Phi(\ox)\ou\ra\le-\nabla\ph_0(\ox)\ou=\la a_j,\nabla\Phi(\ox)\ou\ra\;\mbox{ for }\;j\in J_1.
$$
Employing this along with (\ref{eq35}) tells us that $\nabla\Phi(\ox)\ou\in{\cal K}(\oz,\ov)$ with $\ou\ne 0$. Since we have $\ov=\ov_1+\ov_2$, it follows from (\ref{eq31}) and (\ref{eq34}) that
$$
\big(\ph_0(x_k)-\ph_0(\ox)\big)+\la\ov,z_k-\oz\ra\le 0.
$$
The latter implies by the Taylor expansion together with $\ov\in\Lambda_{\small\mbox{com}}(\ox)$ that
$$
\la\nabla_{xx}^2L(\ox,\ov)\ou,\ou\ra\le 0\;\mbox{ with }\;\ou\ne 0,
$$
which contradicts (\ref{sosc2}) and thus verifies the strict local optimality of $\ox$ in \eqref{comp}.

It remains to justify the noncriticality of $\ov\in\Lambda_{\small\mbox{com}}(\ox)$ in \eqref{comp}. It follows from the discussion before the formulation of the theorem that it suffices to show that condition \eqref{sosc} ensures the noncriticality of the corresponding vector $\ov\in\Lambda(\ox)$ in \eqref{VS}. Assuming the contrary and applying Theorem~\ref{cric002}, we conclude that there is a pair $(\xi,\eta)\in\R^n\times\R^m$ with $\xi\ne 0$ satisfying all the conditions in \eqref{crc6}. This gives us the relationships
$$
0=\la\nabla_x\Psi(\ox,\ov)\xi+\nabla\Phi(\ox)^*\eta,\xi\ra=\la\nabla_x\Psi(\ox,\ov)\xi,\xi\ra+\la\eta,\nabla\Phi(\ox),\xi\ra=\la\nabla_x\Psi(\ox,\ov)\xi,\xi\ra,
$$
which contradict \eqref{sosc} with $u=\xi$ and therefore complete the proof of the theorem. $\h$\vspace*{-0.1in}

\begin{Remark}{\bf(on second-order sufficient conditions and the existence of Lagrange multipliers in composite optimization.)}\label{socs4} {\rm The following discussions on the assumptions and conclusions of Theorem~\ref{sosc3} are useful.

{\bf(i)} Another type of second-order sufficient conditions for strict minimizers of a general extended-real-valued function $\ph\colon\R^n\to\oR$ is derived in \cite[Theorem~13.24]{rw} in terms of the (directional) {\em second subderivative} $\d^2\ph$ taken from \cite[Definition~13.3]{rw}. Applications of this result to  structural optimization problems like the one \eqref{comp} of our consideration require second-order calculus rule for $d^2\ph$. In particular, the chain rule for $d^2(\th\circ\Phi)$ from \cite[Theorem~13.14]{rw} can be implemented in \eqref{comp} by using the calculation of the critical cone for CPWL functions, which is done above in Proposition~\ref{ccid}. On the other hand, the crucial chain rule of \cite[Theorem~13.14]{rw} is obtained under the {\em basic qualification condition}
\begin{equation}\label{rcq}
\sub^{\infty}\th(\oz)\cap\ker\nabla\Phi(\ox)^*=\{0\}\;\mbox{ with }\;\oz=\Phi(\ox),
\end{equation}
which is actually equivalent to Robinson's constraint qualification (RCQ) (see, e.g., \cite[Definition~2.86]{bs}) for problems with the constraints $\Phi(x)\in\dom\th$ as in \eqref{comp}. Note that \eqref{rcq} ensures that $\Lambda_{\small{\rm com}}(\ox)\ne\emp$  , which is an assumption of Theorem~\ref{sosc3}.

{\bf(ii)} The qualification condition \eqref{rcq} is a major requirement for (fully) {\em amenable} compositions as in \cite[Definition~10.23]{rw}, which is {\em not} imposed in our Theorem~\ref{sosc3}. It is a direct consequence of the Mordukhovich criterion \cite[Theorem~9.40]{rw} to see that \eqref{rcq} is equivalent to the {\em metric regularity} of the set-valued mapping
\begin{equation}\label{metreg}
(x,\al)\mapsto F(x,\al):=\epi\th-\big(\Phi(x),\al\big)
\end{equation}
around $(\ox,\oz,0,0)\in\R^n\times\R\times\R^m\times\R$. However, the assumption $\Lambda_{\small{\rm com}}(\ox)\ne\emp$ of Theorem~\ref{sosc3} is satisfied under less restrictive qualification conditions; in particular, under the {\em metric subregularity} of mapping \eqref{metreg} at $(\ox,\oz,0,0)$ (equivalent to the {\em calmness} of its inverse); see \cite{ho05,io} for more details. This allows us to invoke the (fully) {\em subamenable} \cite{gm16} (vs.\ amenable) property of the constraint set $\Phi(x)\in\dom\th$ in \eqref{comp} to get $\Lambda_{\small{\rm com}}(\ox)\ne\emp$. Note that there are a number of other {\em constraint qualifications} ensuring the latter requirement for particular classes of composite optimization problems, especially, for NLPs; see, e.g., \cite{bs,gm15,is14} and the references therein.}
\end{Remark}\vspace*{-0.3in}

\section{Noncriticality from Full Stability in Composite Optimization}\sce\vspace*{-0.05in}

In this section we consider the two-parametric version of problem (\ref{comp}) given by
\begin{equation}\label{comp2}
\mbox{minimize }\;\ph_0(x,p_2)+\theta\big(\Phi(x,p_2)\big)-\la p_1,x\ra\;\mbox{ subject to }\;x\in\R^n
\end{equation}
with $(p_1,p_2)\in\R^n\times\R^l$. Fix $\gg>0$ and $(\ox,\op_1,\op_1)$ with $\Phi(\ox,\op_2)\in\dom\th$ and then define the parameter-depended optimal value function for (\ref{comp}) by
\begin{eqnarray*}
m_\gg(p_1,p_2):=\inf_{\|x-\ox\|\le\gg}\big\{\ph_0(x,p_2)+\theta(\Phi(x,p_2))-\la p_1,x\ra\big\}
\end{eqnarray*}
and the parameterized set of optimal solutions to (\ref{comp}) by
\begin{eqnarray*}
M_\gg(p_1,p_2):=\mbox{argmin}\big\{\ph_0(x,p_2)+\theta\big(\Phi(x,p_2)\big)-\la p_1,x\ra\big|\;\|x-\ox\|\le\gg\big\}
\end{eqnarray*}
with the convention that argmin:=$\emp$ when the expression under minimization is $\infty$. According to the scheme of \cite{lpr} suggested for general optimization problems with extended-real-valued objectives, we say that $\ox$ is a {\em fully stable} locally optimal solution to problem (\ref{comp2}) if there exist a number $\gg>0$ and neighborhoods $U$ of $\op_1$ and $W$ of $\op_2$ such that the mapping $(p_1,p_2)\mapsto M_\gg(p_1,p_2)$ is single-valued and Lipschitz continuous with $M_\gg(\op_1,\op_2)=\{\ox\}$ and the function $(p_1,p_2)\mapsto m_\gg(p_1,p_2)$ is likewise Lipschitz continuous on $U\times W$.

In what follows we concentrate on the {\em canonically perturbed} case of \eqref{comp2} described by
\begin{equation}\label{comp3}
\mbox{minimize }\;\ph_0(x)+\theta\big(\Phi(x)+p_2\big)-\la p_1,x\ra\;\mbox{ subject to }\;x\in\R^n
\end{equation}
with $(p_1,p_2)\in\R^n\times\R^m$ and suppose that the function $\ph_0$ and the mapping $\Phi$ are ${\cal C}^2$-smooth around the reference points. The next theorem shows that full stability of the given locally optimal solution $\ox$ to \eqref{comp3} with $\th\in CPWL$ rules out the existence of critical multipliers associated with $\ox$. This {\em proves the conjecture} of \cite{m15} for the class of composite optimization problems \eqref{comp} studied in the paper; see Section~1 for more discussions.\vspace*{-0.1in}

\begin{Theorem}{\bf(excluding critical multipliers by full stability).}\label{fucri} Let $\ox$ be a fully stable locally optimal solution to {\rm(\ref{comp3})} with $(\op_1,\op_2)=(0,0)$, and let $\th\in CPWL$. Then the Lagrange multiplier set $\Lm_{\small{\rm com}}(\ox)$ in {\rm (\ref{laset2})} does not include any critical multipliers.
\end{Theorem}\vspace*{-0.05in}
{\bf Proof}. We first verify that the imposed full stability of $\ox$ implies the validity of the qualification condition (\ref{rcq}). To proceed, pick any $\eta\in\sub^{\infty}\th(\oz)\cap\ker\nabla\Phi(\ox)^*$ and deduce from the convexity of $\th$ that $\sub^{\infty}\th(\oz)=N(\oz;\dom\th)$. Select $p_1=\op_1=0$ and $p_2=t\eta$ with $t\downarrow 0$. The property of full stability for $\ox$ allows us to find a Lipschitz constant $\ell\ge 0$ and a unique solution to problem (\ref{comp3}), denoted by $x_{p_1p_2}$, for which the following holds:
$$
\|x_{p_1p_2}-\ox\|\le\ell\|p_2\|=\ell t\|\eta\|.
$$
By $\Phi(x_{p_1p_2})+p_2\in\dom\th$ and $\Phi(x_{p_1p_2})=\Phi(\ox)+\nabla\Phi(\ox)(x_{p_1p_2}-\ox)+o(\|x_{p_1p_2}-\ox\|)$ we have
$$
\begin{array}{lll}
0&\ge&\la\eta,\Phi(x_{p_1p_2})+p_2-\Phi(\ox)\ra\\
&=&\la\eta,\nabla\Phi(\ox)(x_{p_1p_2}-\ox)+o(\|x_{p_1p_2}-\ox\|)+p_2\ra\\
&=&\la\eta,o(\|x_{p_1p_2}-\ox\|)\ra+t\|\eta\|^2,
\end{array}
$$
which tells us that $\eta=0$ and thus justifies the validity of (\ref{rcq}).

Pick now $\ov\in\Lm_{\small\mbox{com}}(\ox)$ and prove that $\ov$ is noncritical. Consider the KKT system for problem (\ref{comp3}) and write it in the following form of the canonically perturbed generalized equation:
\begin{equation}\label{geq4}
\left[\begin{array}{c}
p_1\\p_2\\
\end{array}
\right]
\in\left[\begin{array}{c}
\nabla_xL(x,v)\\-\Phi(x)
\end{array}
\right]+\left[\begin{array}{c}
0\\(\sub\th)^{-1}(v)
\end{array}
\right]
\end{equation}
Denote by $S_{\tiny\mbox{KKT}}\colon\R^n\times\R^m\tto\R^n\times\R^m$ the solution map to \eqref{geq4} defined as
\begin{equation}\label{s-kkt}
S_{\tiny\mbox{KKT}}(p_1,p_2):=\big\{(x,v)\in\R^n\times\R^m\big|\;p_1=\nabla_x L(x,v),\;v\in\partial\th\big(p_2+\Phi(\ox)\big)\big\}.
\end{equation}
By Theorem~\ref{uplip} it suffices to show that there are numbers $\ve>0$ and $\ell\ge 0$ as well as neighborhoods $U$ of $0\in\R^n$ and $W$ of $0\in\R^m$ such that for any $(p_1,p_2)\in U\times W$ and any $(x_{p_1p_2},v_{p_1p_2})\in S_{\tiny\mbox{KKT}}(p_1,p_2)\cap(\B_\ve(\ox)\times\B_\ve(\ov))$  we have estimate (\ref{upper}) with replacing $\Lambda(\ox)$ by $\Lm_{\small\mbox{com}}(\ox)$. To verify it, remember that $\ox$ is a fully stable local minimizer of problem (\ref{comp3}) and then deduce from \cite[Proposition~6.1]{mrs} that there are neighborhoods $\Tilde U\times\Tilde W$ of $(0,0)$ and $\Tilde V$ of $\ox$ for which  the set-valued mapping
$$
(p_1,p_2)\mapsto Q(p_1,p_2):=\big\{x\in\R^n\big|\;p_1\in\nabla\ph_0(x)+\nabla\Phi(x)^*\sub\th\big(\Phi(x)+p_2\big)\big\}
$$
admits a Lipschitzian single-valued graphical localization on $\Tilde U\times\Tilde W\times\Tilde V$, which amounts to saying that there exists a Lipschitzian single-valued mapping $g\colon\Tilde U\times\Tilde W\to\Tilde V$ such that $(\gph Q)\cap(\Tilde U\times\Tilde W\times\Tilde V)=\gph g$. Denote $U:=\Tilde U$, $W:=\Tilde W$ and take $\ve>0$ so small that $\B_\ve(\ox)\subset\Tilde V$. By the Lipschitzian single-valued graphical localization property of $Q$, find a constant $\ell\ge 0$ such that for any $(p_1,p_2)\in U\times W$ and any $(x_{p_1p_2},v_{p_1p_2})\in S_{\tiny\mbox{KKT}}(p_1,p_2)\cap(\B_\ve(\ox)\times\B_\ve(\ov))$ it follows that $x_{p_1p_2}\in Q(p_1,p_2)$, and therefore we arrive at the estimate
$$
\|x_{p_1p_2}-\ox\|=\|x_{p_1p_2}-x_{\op_1\op_2}\|\le\ell\big(\|p_1\|+\|p_2\|\big).
$$
As shown in the proof of Theorem~\ref{uplip}, estimate (\ref{up2}) holds with replacing $\Lm(\ox)$ by $\Lm_{\small\mbox{com}}(\ox)$. Adjusting finally $\ve$ if necessary, we get (\ref{upper}) and complete the proof of the theorem. $\h$\vspace*{0.03in}

Theorem~\ref{fucri} extends to the general case of $\th\in CPWL$ in \eqref{geq4} the result by Izmailov \cite{iz15} obtained for classical nonlinear programs with replacing \eqref{rcq} by the Mangasarian-Fromovitz constraint qualification (MFCQ), which is equivalent to \eqref{rcq} for NLPs. Furthermore, Izmailov \cite[Example~3.2]{iz15} constructed the following NLP example:
\begin{equation}\label{iz-ex}
\mbox{minimize }\;x_1+x^4_2\;\mbox{ subject to }\;-x_1\le 0,\;(x_1-2)^2+x^2_2\le 4,\;x=(x_1,x_2)\in\R^2,
\end{equation}
where the canonical perturbation of {\em only the constraints} while not of the cost function (i.e., when $p_1=0$ in \eqref{comp3}) did not guarantee the noncriticality of Lagrange multipliers associated with the local minimizer $\ox=0$, even under the validity of MFCQ at $\ox$.

Observe finally that the result of Theorem~\ref{fucri} allows us to make a conclusion that ``bad" critical multipliers associated with a given local minimizer of \eqref{comp} will not appear (and hence convergent primal-dual algorithms to find this minimizer exhibit {\em high convergent rates}) while operating entirely with the {\em initial data} of \eqref{comp}. It is due to characterizations of full stability for various subclasses of \eqref{comp} with $\th\in CPWL$ obtained recently in \cite{mn14,mn15,mn16,mnr,mrs,ms15a,ms152}.\vspace*{-0.15in}

\section{Tilt Stability versus Critical Multipliers}\sce\vspace*{-0.05in}

This section concerns another challenging issue that was brought up in \cite{m15} about efficient conditions under which critical multipliers are ruled out by {\em tilt stability} of local minimizers, a weaker property than its full stability counterpart, which corresponds to the canonical perturbation of {\em only the cost function} in \eqref{comp3}, i.e., when $p_2=0$ therein. We consider again the composite optimization framework \eqref{comp} with $\th\in CPWL$ and suppose without loss of generality that $0\in\aff\partial\th(\oz)$, where ``aff" stands for the affine hull of the set. As shown in \cite[Section~3]{ms152}, the latter assumption does not indeed impose any restrictions to our second-order analysis.

It is proved in \cite[Lemma~3.1]{ms152} that for any CPWL function $\th\colon\R^m\to\oR$ there exist a positive number $s\le m$, an $s\times m$-matrix $B$, and a CPWL function $\vartheta\colon\R^s\to\oR$, all constructively built via the initial data of $\th$ in \eqref{eq00}, for which
$$
\th(z)=(\vartheta\circ h)(z)\;\mbox{ with }\;h(z):=Bz\;\mbox{ for all }\;z\;\mbox{ around }\;\oz.
$$
Using this, we say that  $\ox$ is a {\em nondegenerate point} of $\Phi$ from (\ref{comp}) relative to $h(z)=Bz$ if
\begin{equation}\label{fnond}
\nabla\Phi(\ox)\R^n+\ker B=\R^m\;\mbox{ with }\;\oz=\Phi(\ox).
\end{equation}
The reader is referred to \cite{ms152} for more details on \eqref{fnond} and its applications. The reader can find therein that \eqref{fnond} can be equivalently written in the form
\begin{equation}\label{nondeg}
\aff\partial\th(\oz)\cap\ker\nabla\Phi(\ox)^*=\big\{0\big\}.
\end{equation}

Here we employ \eqref{fnond} to establish relationships between critical multipliers and tilt stability of local optimal solutions in composite optimization. Following the line of Poliquin and Rockafellar \cite{pr98}, consider the one-parametric problem
\begin{equation}\label{comp05}
\mbox{minimize }\;\ph_0(x)+\theta\big(\Phi(x)\big)-\la p_1,x\ra\;\mbox{ subject to }\;x\in\R^n
\end{equation}
with $p_1\in\R^n$ and the solution map to it defined by
\begin{eqnarray*}
p_1\mapsto M_\gg(p_1):=\mbox{argmin}\big\{\ph_0(x)+\theta\big(\Phi(x)\big)-\la p_1,x\ra\big|\;\|x-\ox\|\le\gg\big\}.
\end{eqnarray*}
Then $\ox$ is a {\em tilt-stable} local minimizer of \eqref{comp05} if the solution map $M(\cdot)$ is locally single-valued and Lipschitz continuous around $(0,\ox)$ with $M_\gg(0)=\{\ox\}$. Tilt stability is clearly a particular case of full stability. The following result is a consequence of Theorem~\ref{fucri}. \vspace*{-0.1in}

\begin{Corollary}{\bf(excluding critical multipliers by tilt stability under nondegeneracy).}\label{tilt} Let $\ox$ be a tilt-stable locally optimal solution to problem  \eqref{comp05} with $\op_1=0\in\R^n$, let $\th\in CPWL$, and let the nondegeneracy condition \eqref{fnond} be satisfied. Then the Lagrange multiplier set $\Lm_{\small{\rm com}}(\ox)$ from \eqref{laset2} is singleton and the unique Lagrange multiplier in $\Lm_{\small{\rm com}}(\ox)$ is noncritical.
\end{Corollary}\vspace*{-0.05in}
{\bf Proof}. The uniqueness of Lagrange multipliers is a consequence of nondegeneracy \eqref{fnond} and is proved in \cite[Proposition~4.2]{ms152}. To justify the noncriticality of the unique Lagrange multiplier, we invoke \cite[Lemma~6.1]{ms152}, which tells us that under \eqref{fnond} the tilt stability of $\ox$ in problem (\ref{comp05}) is equivalent to the full stability of the canonically perturbed problem (\ref{comp2}) at this point. Employing now Theorem~\ref{fucri} yields the noncriticality of $\ox$. $\h$\vspace*{-0.1in}

\begin{Remark}{\bf(excluding critical multipliers by tilt stability for degenerate NLPs).}\label{tilt-deg} {\rm In the particular case of {\em NLPs} we can conclude that tilt stability of a local minimizer $\ox$ excludes the existence of critical multipliers associated with it under {\em weaker} qualification conditions than the nondegeneracy in Corollary~\ref{tilt}. Indeed, it is shown in \cite[Theorem~4.3]{mn15} that tilt stability of $\ox$ is equivalent to the so-called {\em uniform second-order sufficient condition} (USOSC) under the simultaneous validity of MFCQ and the constant rank constraint qualification (CRCQ). Since USOSC yields the classical second-order sufficient condition in NLPs, which in turns clearly excludes criticality of all the multipliers associated with $\ox$, we get that the {\em combination of MFCQ and CRCQ},  which is strictly weaker than nondegeneracy, ensures the noncriticality of all the corresponding multiplies. Furthermore, the second-order characterization of tilt stability for NLPs obtained in \cite[Theorem~7.6]{gm15} via only the {\em extreme} multipliers in {\em critical directions} at $\ox$ allows us exclude criticality of {\em all} the multipliers associated with this local minimizer.}
\end{Remark}\vspace*{-0.1in}

As we see from the very construction of \eqref{comp05}, tilt stability reflects only the cost function perturbation without any perturbation of the constraints. The following two examples show that it may not be possible to rule out critical multipliers, even under plausible constraint qualifications that are weaker than the simultaneous validity of MFCQ and CRCQ. These examples are complementary to the one by Izmailov \eqref{iz-ex} showing that only the constraint perturbations are not sufficient for noncriticality. Both examples below are slight modifications of those in \cite[Examples~8.3, 8.4]{gm15}, which were constructed there for different purposes.\vspace*{-0.1in}

\begin{Example}{\bf(presence of critical multipliers for tilt-stable minimizers of NPLs under SOSCMS).}\label{ex01}
{\rm Consider the three-dimensional nonlinear program:
\begin{equation*}
\left\{\begin{array}{ll}
{\rm{minimize}}&\ph_0(x):=\disp-x_1+\frac{5}{2}x_2^2+x_3^2,\quad x=(x_1,x_2,x_3)\in\R^3,\\
{\rm{subject\;to}}&\ph_1(x):=\disp x_1-\frac{1}{2}x_2^2\le 0,\;\;\ph_2(x):=x_1-\frac{1}{2}x_3^2\le 0,\\
&\ph_3(x):=-x_1-\disp\frac{1}{2}x_2^2-\frac{1}{2}x_3^2\le 0.
\end{array}\right.
\end{equation*}
As follows from \cite[Theorem~7.6]{gm15}, the local minimizer $\ox=0$ is tilt-stable for this problem due to the validity of the {\em second-order sufficient condition for metric subregularity} (SOSCMS) by Gfrerer \cite{g}, although both MFCQ and CRCQ fail at $\ox$. Direct calculations show that the set of Lagrange multipliers associated with $\ox$ is
$$
\Lambda(\ox)=\big\{v=(v_1,v_2,v_3)\in\R_+^3\big|\;v_1+v_2-v_3=0\big\}
$$
and that the Lagrangian Hessian at $\ox$ and the multiplier $\ov=(3,0,2)$ is
$$
\nabla_{xx}^2L(\ox,\ov)=\nabla^2\ph_0(\ox)+\ov_1\nabla^2\ph_1(\ox)+\ov_2\nabla^2\ph_2(\ox)+\ov_3\nabla^2\ph_3(\ox)=0.
$$
Then we observe that the pair $(\xi,\eta)$ with $\xi=(0,1,1)$ and $\eta=(0,0,0)$ satisfies all the conditions in (\ref{forex1}) with $\Psi=\nabla_x L$. This confirms by Theorem~\ref{cric002} and Example~\ref{nlpineq} that $\ov$ is a critical multiplier at $\ox$ for the nonlinear program (\ref{ex01}).}
\end{Example}\vspace*{-0.1in}

The next example shows that the {\em MFCQ alone} may not rule out the existence of critical multipliers in three-dimensional NLPs.\vspace*{-0.1in}

\begin{Example}{\bf(presence of critical multipliers for tilt-stable minimizers of NLPs under MFCQ).}\label{ex2}
{\rm Consider the following three-dimensional nonlinear program:
\begin{equation*}
\left\{\begin{array}{ll}
{\rm{minimize}}&\disp\ph_0(x):=-x_1+\frac{1}{2}x_2^2,\quad x=(x_1,x_2,x_3)\in\R^3,\\
{\rm{subject\;to}}& \ph_1(x):=x_1+x_3^2\le 4,\quad\ph_2(x):=x_1\le 0.
\end{array}\right.
\end{equation*}
Arguing as in \cite[Example~8.4]{gm15} by using \cite[Theorem~6.1]{gm15} tells us that $\ox=(0,0,0)$ is a tilt-stable minimizer for this problem. We easily check that MFCQ is fulfilled at $\ox$ and determine that the set of Lagrange multipliers at this minimizer is
$$
\Lambda(\ox)=\big\{v=(v_1,v_2)\in\R_+^2\big|\;v_1+v_2=1\big\}.
$$
Take $\ov=(\ov_1,\ov_2)=(0,1)\in\Lambda(\ox)$ and get $I_+(\ov)=\{2\}$ and $I_0(\ov)=\{1\}$ in the notation of Example~\ref{nlpineq}. Further, we calculate the Lagrangian Hessian at $(\ox,\ov)$ by
$$
\nabla_{xx}^2L(\ox,\ov)=\nabla^2\ph_0(\ox)+\ov_1\nabla^2\ph_1(\ox)+\ov_2\nabla^2\ph_2(\ox)=\left(\begin{array}{ccc}
0&0&0\\0&1&0\\0&0&0
\end{array}\right)
$$}
\end{Example}
and deduce that the pair $(\xi,\eta)$ with $\xi=(0,0,1)$ and $\eta=(0,0,0)$ satisfies all the conditions in (\ref{forex1}) with $\Psi=\nabla_x L$. This confirms by Theorem~\ref{cric002} the criticality of $\ov$ at $\ox$.\vspace*{-0.15in}

\section{Robust Isolated Calmness via Noncritical Multipliers}\sce\vspace*{-0.05in}

This section concerns some stability properties of set-valued mappings, which were first designated by Robinson \cite{rob} under the ``upper Lipschitzian" name and then has been widely spread in variational analysis under the name of ``calmness" or ``metric subregularity" of the inverse; see, e.g., \cite{rw}. A further specified notion of this type scattered in variational analysis under different names is finally formulated as follows \cite{dr}: A mapping $F\colon\R^n\tto\R^m$ is said to be {\em isolatedly calm} at $(\ox,\oy)\in\gph F$ if there exist a constant $\ell\ge0$ and neighborhoods $U$ of $\ox$ and $V$ of $\oy$ such that
\begin{equation}\label{calm}
F(x)\cap V\subset\{\oy\}+\ell\|x-\ox\|\B\;\mbox{ for all }\;x\in U.
\end{equation}
The isolated calmness property \eqref{calm} admits the following characterization via the graphical derivative \eqref{gder}, the necessity part of which was obtained in \cite[Proposition~2.1]{kr} while the sufficiency was proved later in \cite[Proposition~4.1]{l}:
\begin{equation}\label{calm-cr}
DF(\ox,\oy)(0)=\{0\}.
\end{equation}
Obstacles with applications of \eqref{calm} arise from {\em nonrobustness} of the graphical derivative, as well as of the calmness and isolated calmness properties themselves. Nevertheless, recent results on the calculation of the graphical derivative for some particular mappings describing solution maps to certain kinds of generalized equations have generated by \eqref{calm-cr} efficient conditions for isolated calmness of perturbed variational systems associated with constraints of the type $g(x)\in\Th$ for smooth $g$ under various  qualification conditions and assumptions on $\Th$; see \cite{ch,dr,go,mor1,mor2}.

Quite recently \cite{dsz}, a robust version of \eqref{calm} with the additional requirement that $F(x)\cap V\ne\emp$ for all $x\in U$ has been labeled as the {\em robust isolated calmness} of $F$ at $(\ox,\oy)$. Note that this property was actually employed earlier in particular settings under different names or without naming it at all; see \cite{bo94,dr97,kk,p}. If the set-valued mapping $F$ is lower semicontinuous at $(\ox,\oy)$ in the standard topological sense, then isolated calmness implies its robust counterpart. However, it does not hold in general as shown, e.g., in \cite[Example~6.4]{mor2}.

It is worth mentioning that the usage of robust isolated calmness in numerical optimization has been recognized in the literature starting with 1990s. In particular, the sharpest result for the sequential quadratic programming (SQP) method for solving NLPs, obtained by Bonnans \cite{bo94}, imposes the strict Mangasarian-Fromovitz constraint qualification together with the conventional second-order sufficient condition for NLPs. As later proved by Dontchev and Rockafellar \cite[Theorem~2.6]{dr97}, the simultaneous validity of these conditions characterizes the robust isolated calmness of solutions maps of canonically perturbed KKT systems in NLPs. Recently this result has been extended by Ding et al. \cite[Theorem~24]{dsz} to some nonpolyhedral problems of constrained optimization under the so-called strict Robinson constraint qualification.

The main goal of this section is to establish close relationships between isolated calmness and its robust counterpart for solutions maps to canonically perturbed variational systems \eqref{geq2} from one side and {\em noncritical multipliers} from the other, which do not seem to be explicitly recognized in the literature even for NLPs. We obtain such relationships in the general variational setting of \eqref{VS} and then efficiently specify and strengthen them in the case of KKT systems in composite optimization \eqref{comp} with $\th\in CPWL$. As a by-product of these developments, we offer a {\em new viewpoint} on the study of isolated calmness and its robust version in constrained optimization that is essentially different from those developed in \cite{dsz,dr97}.

First we study relationships between noncriticality of multipliers from Definition~\ref{crit} and isolated calmness of the solution map \eqref{mapS} to the canonically perturbed system \eqref{geq2} for general variational systems \eqref{VS} with arbitrary (proper) functions $\th\colon\R^m\to\oR$ and also for $\th\in CPWL$.\vspace*{-0.1in}

\begin{Theorem}{\bf(relationships between noncriticality and isolated calmness for general variational systems).}\label{crit-calm} The noncriticality of the multiplier $\ov\in\Lambda(\ox)$ for \eqref{VS} in the framework of Definition~{\rm\ref{crit}} is equivalent to the implication
\begin{equation}\label{crc5}
\left\{\begin{array}{ll}
\nabla_{x}\Psi(\ox,\ov)\xi+\nabla\Phi(\ox)^*\eta=0,\\
\eta\in\big(D\sub\th\big)(\oz,\ov)\big(\nabla\Phi(\ox)\xi\big)
\end{array}\right.
\Longrightarrow\xi=0,
\end{equation}
while the isolated calmness at $((0,0),(\ox,\ov))\in\gph S$ of the solution map \eqref{mapS} to the canonically perturbed system \eqref{geq2} amounts to the stronger implication
\begin{equation}\label{bn02}
\left\{\begin{array}{ll}
\nabla_{x}\Psi(\ox,\ov)\xi+\nabla\Phi(\ox)^*\eta=0,\\
\eta\in\big(D\sub\th\big)(\oz,\ov)\big(\nabla\Phi(\ox)\xi\big)
\end{array}\right.
\Longrightarrow(\xi,\eta)=(0,0)
\end{equation}
with $S(0,0)\cap V=\{(\ox,\ov)\}$ for some neighborhood $V$ of $(\ox,\ov)$. If furthermore $\th\in CPWL$, then the noncriticality of any $\ov\in\Lm(\ox)$ yields the existence of a neighborhood $V$ of $(\ox,\ov)$ for which
\begin{equation}\label{calm-lagr}
S(0,0)\cap V=\big[\{\ox\}\times\Lambda(\ox)\big]\cap V.
\end{equation}
\end{Theorem}\vspace*{-0.05in}
{\bf Proof.} It follows from the the conditions in \eqref{crc} and from the structure of the mapping $G$ in \eqref{F-KKT} that the noncriticality of $\ov\in\Lm(\ox)$ for \eqref{VS} can be expressed in the form
\begin{equation*}
(0,0)\in DG\big((\ox,\ov),(0,0)\big)\big(\xi,\eta)\Longrightarrow\xi=0\;\mbox{ for }\;(\xi,\eta)\in\R^n\times\R^m.
\end{equation*}
Since $G$ is represented as $G(x,v)=g(x,v)+Q(v)$ with a smooth mapping $g\colon\R^n\times\R^m\to\R^n\times\R^m$, we easily deduce from definition \eqref{gder} of the graphical derivative that
\begin{equation}\label{bn01}
\begin{array}{lll}
DG\big((\ox,\ov),(0,0)\big)(\xi,\eta)&=&\left[\begin{array}{cc}
\nabla_{x}\Psi(\ox,\ov)&\nabla\Phi(\ox)^*\\
-\nabla\Phi(\ox)&0\\
\end{array}
\right]\left[\begin{array}{c}
\xi\\\eta\\
\end{array}\right]
+\left[\begin{array}{c}
0\\D(\sub\th)^{-1}(\ov,\oz)(\eta)\\
\end{array}
\right]\\\\
&=&\left[\begin{array}{c}
\nabla_{x}\Psi(\ox,\ov)\xi+\nabla\Phi(\ox)^*\eta,\\
-\nabla\Phi(\ox)\xi+D(\sub\th)^{-1}(\ov,\oz)(\eta)
\end{array}\right]
\end{array}
\end{equation}
and therefore arrive in this way to the noncriticality characterization \eqref{crc5}.

Apply further the isolated calmness criterion \eqref{calm-cr} to the solution map $F:=S$ from \eqref{mapS} at the point $((0,0),(\ox,\ov))\in\gph S$. Taking into account that the equivalence
\begin{equation*}
(\xi,\eta)\in DS\big((0,0),(\ox,\ov)\big)(w_1,w_2)\Longleftrightarrow(w_1,w_2)\in DG\big((\ox,\ov),(0,0)\big)(\xi,\eta)
\end{equation*}
is valid for any $(\xi,\eta),(w_1,w_2)\in\R^n\times\R^m$ and using the graphical derivative representation \eqref{bn01}, we conclude that the isolated calmness of $S$ at $((0,0),(\ox,\ov))$ reduces to the fulfillment of \eqref{bn02}. The fact on $S(0,0)\cap V=\{(\ox,\ov)\}$ follows directly from the isolated calmness property of $S$.

It remains to verify \eqref{calm-lagr} in the case where $\th\in CPWL$. Using the characterization from Theorem~\ref{uplip} of noncriticality for any $\ov\in\Lambda(\ox)$ in this case, we find numbers $\ve>0$, $\ell\ge 0$ and a neighborhood $U$ of $(0,0)\in\R^n\times\R^m$ such that the estimate
\begin{equation*}
\|x-\ox\|+\dist\big(v;\Lm(\ox)\big)\le\ell\big(\|p_1\|+\|p_2\|\big)
\end{equation*}
holds whenever $(p_1,p_2)\in U$ and any $(x,v)\in S(p_1,p_2)\cap\B_\ve(\ox,\ov)$. Let us justify \eqref{calm-lagr}  by showing that $S(0,0)\cap\B_\ve(\ox,\ov)=[\{\ox\}\times\Lambda(\ox)]\cap \B_\ve(\ox,\ov)$. Indeed, the inclusion
\begin{equation*}
\big[\{\ox\}\times\Lambda(\ox)\big]\cap\B_\ve(\ox,\ov)\subset S(0,0)\cap\B_\ve(\ox,\ov)
\end{equation*}
is a direct consequence of the feasibility of $(\ox,\ov)$ for the variational system \eqref{VS} and the definitions of $S$ and $\Lambda(\ox)$ in \eqref{mapS} and \eqref{laset}, respectively. To get the opposite inclusion, pick any pair $(x,v)\in S(0,0)\cap\B_\ve(\ox,\ov)$ and deduce from the estimate above that
\begin{equation*}
\|x-\ox\|+\dist\big(v;\Lm(\ox)\big)\le\ell\big(\|0\|+\|0\|\big)=0,
\end{equation*}
which tells us that $\|x-\ox\|+\dist(v;\Lm(\ox))=0$. Thus we arrive at the claimed conditions $x=\ox$ and $v\in\Lambda(\ox)$ and complete the proof of theorem.  $\h$\vspace*{-0.05in}

\begin{Remark}{\bf (relationships between calmness, isolated calmness, and noncriticality).}\label{rem07}
{\rm It is worth highlighting the differences between the calmness, isolated calmness, and its semi-isolated version for the mapping $S$ at $((0,0),(\ox,\ov))$ and the noncriticality of the multiplier $\ov\in\Lambda(\ox)$ for \eqref{VS}. The calmness property of $S$ at $((0,0),(\ox,\ov))$ reads as the existence of $\ell\ge 0$ and neighborhoods $U$ of $(0,0)$ and $V$ of $(\ox,\ov)$ so that for any $(p_1,p_2)\in U$ the inclusion
\begin{equation}\label{ca}
S(p_1,p_2)\cap V\subset S(0,0)+\ell\big(\|p_1\|+\|p_2\|\big)\B
\end{equation}
holds while the noncriticality of $\ov\in\Lambda(\ox)$ for \eqref{VS} is equivalent to the existence of $\ell\ge 0$ and neighborhoods $U$ of $(0,0)$ and $V$ of $(\ox,\ov)$ so that for any $(p_1,p_2)\in U$ we have
\begin{equation}\label{ca2}
S(p_1,p_2)\cap V\subset\{\ox\}\times\Lambda(\ox)+\ell\big(\|p_1\|+\|p_2\|\big)\B
\end{equation}
due to estimate \eqref{upper} in Theorem~\ref{uplip}. Finally, the isolated calmness of $S$ at $((0,0),(\ox,\ov))$ amounts to the existence of $\ell\ge 0$ and neighborhoods $U$ of $(0,0)$ and $V$ of $(\ox,\ov)$ such that
\begin{equation}\label{ca3}
S(p_1,p_2)\cap V\subset\big\{(\ox,\ov)\big\}+\ell\big(\|p_1\|+\|p_2\|\big)\B\;\mbox{ whenever }\;(p_1,p_2)\in U.
\end{equation}
Comparing (\ref{ca})--(\ref{ca3}) brings us to the following implications:
$$
\mbox{isolated calmness}\quad\Longrightarrow\mbox{noncriticality}\quad\Longrightarrow\mbox{ calmness},
$$
which allows us to label property \eqref{upper} equivalent to \eqref{ca2} as ``semi-isolated calmness'' to distinguish it from both isolated calmness and calmness properties for the mapping S. Observe to this end that estimate (\ref{upper}) and its equivalent form \eqref{ca2} can be interpreted as the partial isolated calmness of $S$ with respect to $x$ while being reduced to the full isolated calmness of $S$ when the set of Lagrange multipliers at $\ox$ is the singleton $\Lm(\ox)=\{\ov\}$.}
\end{Remark}\vspace*{-0.05in}

We show below that the results of Theorem~\ref{crit-calm} can be significantly strengthen for the case of KKT systems \eqref{geq4} associated with canonically perturbed composite optimization problems of type \eqref{comp3} where $\th\in CPWL$.  In this case the corresponding solution map $S_{\tiny\mbox{KKT}}\colon(p_1,p_2)\mapsto(x,v)$ is given in \eqref{s-kkt} and the set of Lagrange multipliers $\Lambda_{\small{\rm com}}(\ox)$ is taken from \eqref{laset2}.

To proceed in this direction, we first present a second-order {\em necessary} condition for local optimality in the composite optimization problem \eqref{comp} used in what follows.\vspace*{-0.1in}

\begin{Proposition}{\bf(second-order necessary optimality condition for composite problems).}\label{suff2} Let $\ox$ be a feasible solution to \eqref{comp} with
$\th\in CPWL$, and let the qualification condition {\rm(\ref{rcq})} be satisfied. If $\ox$ is a locally optimal solution to {\rm(\ref{comp})}, then $\Lambda_{\small{\rm com}}(\ox)\ne\emp$ and we have the following second-order optimality condition:
\begin{equation}\label{re30}
\max_{\ov\in\Lambda_{\small{\rm com}}(\ox)}\la\nabla_{xx}^2L(\ox,\ov)u,u\ra\ge 0\;\mbox{ for all }\;0\ne u\in\R^n\;\mbox{ with }\;\nabla\Phi(\ox)u\in{\cal K}(\oz,\ov),
\end{equation}
where the critical cone ${\cal K}(\oz,\ov)$ with $\oz=\Phi(\ox)$ is defined in \eqref{cc} and is calculated in \eqref{cc2} via the given data of the CPWL function $\th$.
\end{Proposition}\vspace*{-0.05in}
{\bf Proof.} As discussed, the constraint qualification \eqref{rcq} yields $\Lambda_{\small{\rm com}}(\ox)\ne\emp$; see Remark~\ref{socs4}. Then we apply \cite[Theorem~13.24]{rw} and proceed similarly to \cite[Example~13.25]{rw} that deals with the constraint $g(x)\in\Th$ described by a ${\cal C}^2$-smooth mapping $g$ and a convex polyhedron $\Th$. In our case we use the critical cone ${\cal K}(\oz,\ov)$ from \eqref{cc}, which allows us to arrive in this way at the claimed second-order necessary optimality condition \eqref{re30}. $\h$\vspace*{0.05in}

Next we derive a useful statement of its own interest revealing that the basic qualification condition \eqref{rcq} must be satisfied for any (proper) {\em convex} function $\th\colon\R^m\to\oR$ provided that the sets of Lagrange multipliers \eqref{laset2} to \eqref{comp} at $\ox$ is a singleton.\vspace*{-0.1in}

\begin{Proposition}{\bf(validity of the basic qualification condition).}\label{bqc} Let $\Phi\colon\R^n\to\R^m$ be differentiable at $\ox$, and let $\th\colon\R^m\to\oR$ be convex and finite at $\oz=\Phi(\ox)$. If $\Lambda_{\small{\rm com}}(\ox)=\{\ov\}$ for \eqref{laset2}, then the basic qualification condition \eqref{rcq} is satisfied.
\end{Proposition}\vspace*{-0.05in}
{\bf Proof.} Suppose on the contrary that \eqref{rcq} fails and find a singular subgradient $\Tilde v\in\partial^\infty\th(\oz)$ such that $\nabla\Phi(\ox)^*\Tilde v=0$ while $\Tilde v\ne 0$. Define the vector $\Hat v=\ov+\Tilde v\ne\ov$ and show that $\Hat v\in\partial^\infty\th(\oz)$. Indeed, it immediately follows from the construction of $\Hat v$ that $\nabla_x L(\ox,\Hat v)=0$ for the Lagrangian \eqref{lagr}. Since $\ov\in\partial\th(\oz)$ and $\th$ is convex, we have
\begin{equation*}
\la\Hat v,z-\oz\ra=\la\ov,z-\oz\ra+\la\Tilde v,z-\oz\ra\le\th(z)-\th(\oz)+\la\Tilde v,z-\oz\ra\;\mbox{ for all }\;z\in\dom\th.
\end{equation*}
On the other hand, $\la\Tilde v,z-\oz\ra\le 0$ whenever $z\in\dom\th$ due to the aforementioned singular subdifferential representation $\partial^\infty\th(\oz)=N(\oz;\dom\th)$ for convex functions and normal cone construction in convex analysis. This shows that $\Hat v\in\partial\th(\oz)$ and hence $\Hat v\in\Lambda_{\small{\rm com}}(\ox))$ by \eqref{laset2}, which contradicts the assumption on $\Lambda_{\small{\rm com}}(\ox)=\{\ov\}$ and thus verifies that \eqref{rcq} holds. $\h$\vspace*{0.03in}

Now we are ready to establish the major result of this section showing that the isolated calmness of the solution map $S_{\tiny\mbox{KKT}}$ at $((0,0),(\ox,\ov))$ associated with a local minimizer $\ox$ is actually equivalent to its {\em robust} isolated calmness and that both these calmness properties reduce to the {\em noncriticality} to the unique multiplier $\ov$. Furthermore, all these properties are characterized by the {\em second-order sufficient condition} (SOSC) in \eqref{comp} defined in \eqref{sosc2} and justified in Theorem~\ref{sosc3} for the strict optimality of $\ox$ in composite optimization.\vspace*{-0.1in}

\begin{Theorem}{\bf(characterization of robust isolated calmness for KKT systems of composite optimization).}\label{calm2} Let $\ox$ be a feasible solution to the unperturbed problem \eqref{comp}, and let $\th\in CPWL$. Then the following assertions are equivalent:

{\bf(i)} The solution map $S_{\tiny{\rm KKT}}$ in \eqref{s-kkt} is robustly isolatedly calm at the point $((0,0),(\ox,\ov))\in\R^{n+m}\times\R^{n+m}$ and $\ox$ is a locally optimal solution to {\rm(\ref{comp})}.

{\bf(ii)} SOSC {\rm(\ref{sosc2})} holds and $\Lm_{\small{\rm com}}(\ox)=\{\ov\}$ for the set of Lagrange multipliers \eqref{laset2}.

{\bf(iii)} $\Lm_{\small{\rm com}}(\ox)=\{\ov\}$, $\ox$ is a locally optimal solution to \eqref{comp}, and $\ov$ is a noncritical multipliers for \eqref{VS} with $\Psi=\nabla_xL$ associated with the solution $\ox$.

{\bf(iv)} $S_{\tiny{\rm KKT}}$ is isolatedly calm at $((0,0),(\ox,\ov))$ and $\ox$ is a locally optimal solution to {\rm(\ref{comp})}.
\end{Theorem}\vspace*{-0.05in}
{\bf Proof.} We begin with verifying (ii)$\Longrightarrow$(iii). Having (ii) and employing Theorem~\ref{sosc3} tell us that $\ox$ is a strict local minimizer of \eqref{comp}. As indicated in Section~4, Theorem~\ref{cric002} applied to for \eqref{comp} ensures that SOSC \eqref{sosc2} yields the noncriticality of $\ov\in\Lm_{\small{\rm com}}(\ox)$, and hence we arrive at (iii).

Suppose next that (iii) holds and then verify (i). It follows from Theorem~\ref{uplip} that there are numbers $\ve>0$ and $\ell\geq0$ together with neighborhoods $U_1$ of $\bar p_1=0\in\R^n$ and $U_2$ of $\bar p_2=0\in\R^m$ so that for any $(p_1,p_2)\in U_1\times U_2$ and $(x_{p_1p_2},v_{p_1p_2})\in S_{\tiny{\rm KKT}}(p_1,p_2)\cap\B_\ve(\ox,\ov)$ we get
\begin{equation}\label{re15}
\|x_{p_1p_2}-\ox\|+\dist\big(v_{p_1p_2};\Lm_{\small{\rm com}}(\ox)\big)\le\ell\big(\|p_1\|+\|p_2\|\big).
\end{equation}
Combining \eqref{re15} with $\Lm_{\small{\rm com}}(\ox)=\{\ov\}$ gives us neighborhoods $V$ of $\ox$ and $W$ of $\ov$ for which
\begin{equation*}
\|x-\ox\|+\|v-\ov\|\le\ell\big(\|p_1\|+\|p_2\|\big)\;\mbox{ if }\;(x,v)\in S_{\tiny{\rm KKT}}(p_1,p_2)\cap(V\times W),\;(p_1,p_2)\in U_1\times U_2.
\end{equation*}
This shows that the solution map $S_{\small{\rm KKT}}$ is isolatedly calm at $((0,0),(\ox,\ov))$. To get (i), it remains to verify that $S_{\small{\rm KKT}}$ is {\em robustly} isolatedly calm at the point.

We proceed by considering the set-valued mapping $H\colon\R^m\tto\R^n$ defined by
\begin{equation*}
H(p):=\big\{x\in\R^n\big|\;\Phi(x)+p\in\dom\th\big\},\quad p\in\R^m.
\end{equation*}
It follows from Proposition~\ref{bqc} that the qualification condition \eqref{rcq} holds by the assumption $\Lm_{\small{\rm com}}(\ox)=\{\ov\}$ in (iii).
Then we can deduce from \cite[Theorem~4.37(ii)]{m06} applied to the mapping $H$ at $(0,\ox)\in\R^m\times\R^n$ that there are numbers $r>0$ and $\ell\ge 0$
such that
\begin{equation}\label{feas2}
H(p)\cap\B_r(\ox)\subset H(p')+\ell\|p-p'\|\B\;\mbox{ for all }\;p,p'\in r\B,
\end{equation}
where $r>0$ is chosen so small that $\B_r(\ox)\subset V$. Consider now the optimization problem:
\begin{equation}\label{comp4}
\mbox{minimize }\;\ph_0(x)+\theta(\Phi(x)+p_2)-\la p_1,x\ra\;\mbox{ subject to }\;x\in\B_r(\ox)\cap H(p_2),
\end{equation}
which clearly admits an optimal solution $x_{p_1p_2}$ for any pair $(p_1,p_2)\in U_1\times U_2$.\\\vspace*{0.05in}
{\bf Claim~1:} {\em There exists $\epsilon>0$ with $\B_\epsilon(0,0)=\ve\B\subset U_1\times U_2$ such that
\begin{equation*}
x_{p_1p_2}\in\int\B_r(\ox)\;\mbox{ for any }\;(p_1,p_2)\in\ve\B.
\end{equation*}}
Indeed, assuming the contrary gives us a sequence $(p_{1k},p_{2k})\to(0,0)$ and a sequence of optimal solutions $x_{p_{1k}p_{2k}}$ to (\ref{comp4}) with
$\|x_{p_{1k}p_{2k}}\|=r$. Considering a subsequence of $\{x_{p_{1k}p_{2k}}\}$ if necessary, suppose that $x_{p_{1k}p_{2k}}\to\tilde{x}$ for some $\Tilde x\ne\ox$ with $\|\tilde{x}\|=r$. The optimality of $x_{p_{1k}p_{2k}}$ in (\ref{comp4}) yields
\begin{equation}\label{re30a}
\ph_0(x_{p_{1k}p_{2k}})+\theta\big(\Phi(x_{p_{1k}p_{2k}})+p_{2k}\big)-\la p_{1k},x_{p_{1k}p_{2k}}\ra\le\ph_0(x)+\theta\big(\Phi(x)+p_{2k}\big)-\la p_{1k},x\ra
\end{equation}
for any $x\in\B_r(\ox)\cap H(p_{2k})$. Let us now show that
\begin{equation}\label{re33}
\ph_0(\tilde{x})+\theta\big(\Phi(\tilde{x})\big)\le\ph_0(x)+\theta(\Phi(x)\big)\;\;\mbox{whenever}\;\;x\in\B_{\frac{r}{2}}(\ox)\cap H(0),
\end{equation}
which contradicts the strict local optimality of $\ox$ for the unperturbed problem (\ref{comp}). To verify \eqref{re33}, pick $x\in\B_{\frac{r}{2}}(\ox)\cap H(0)$ and take $k\in\N$ so large that $p_{2k}\in\al\B$ with $\al<\min\{\frac{r}{2\ell},r\}$. By (\ref{feas2}) we find $x=x'+\ell\|p_{2k}\|b$ with some $x'\in H(p_{2k})$ and $b\in\B$ for which
\begin{equation*}
\|x'-\ox\|\le\|x-\ox\|+\ell\|p_{2k}\|\le\frac{r}{2}+\ell\frac{r}{2\ell}=r.
\end{equation*}
This implies that $x'\in\B_r(\ox)\cap H(p_{2k})$. Substituting $x'$ into (\ref{re30a}) gives us the estimate
\begin{equation*}
\begin{array}{lll}
\ph_0(x_{p_{1k}p_{2k}})+\theta\big(\Phi(x_{p_{1k}p_{2k}})+p_{2k}\big)-\la p_{1k},x_{p_{1k}p_{2k}}\ra &\le\ph_0\big(x-\ell\|p_{2k}\|b\big)\\
&+\theta\big(\Phi(x-\ell\|p_{2k}\|b)+p_{2k}\big)-\la p_{1k},x-\ell\|p_{2k}\|b\ra,
\end{array}
\end{equation*}
which yields (\ref{re33}) by passing to the limit as $k\to\infty$ and thus justifies this claim.\vspace*{0.03in}

To continue the verification of (i), we deduce from Claim~1 that $\Lm_{\small{\rm com}}(x_{p_1p_2})\ne\emp$ for all $(p_1,p_2)\in\ve\B$ when $\ve$ is sufficiently small. This follows from the validity of $\Lm_{\small{\rm com}}(\ox)\ne\emp$ under the qualification condition \eqref{rcq} and its robustness with respect to perturbations of the initial point. Letting $v_{p_1p_2}\in\Lm_{\small{\rm com}}(x_{p_1p_2})$ and arguing as in the proof of (\ref{up2}) via the Hoffman Lemma tell us, when $\epsilon$ is small, that $(x_{p_1p_2},v_{p_1p_2})\in S_{\tiny{\rm KKT}}(p_1,p_2)\cap(U\times W)$ for any $(p_1,p_2)\in\ve\B$, which justifies the robust isolated calmness in (i) and thus completes the proof of (iii)$\Longrightarrow$(i).\vspace*{0.03in}

Let us next prove (iii)$\Longrightarrow$(ii). It follows from the property $\Lm_{\small{\rm com}}(\ox)=\{\ov\}$ in (iii) that the qualification condition \eqref{rcq} holds by Proposition~\ref{bqc}. Since $\ox$ in (iii) is a local minimizer for \eqref{comp}, we get from the second-order necessary optimality condition of Proposition~\ref{suff2} that
\begin{equation}\label{re17}
\la\nabla_{xx}^2L(\ox,\ov)u,u\ra\ge 0\;\mbox{ for all }\;0\ne u\in\R^n\;\mbox{ with }\;\nabla\Phi(\ox)u\in{\cal K}(\oz,\ov),
\end{equation}
where the critical cone ${\cal K}(\oz,\ov)$ is taken from (\ref{cc2}). To obtain the remaining SOSC in (ii), let us check that the noncriticality of $\ov$ in (iii) ensures that the inequality in \eqref{re17} is strict for $u\ne 0$.\\\vspace*{0.05in}
{\bf Claim~2.} {\em If there is $\ou\ne 0$ satisfying $\nabla\Phi(\ox)\ou\in{\cal K}(\oz,\ov)$ and
$\la\nabla_{xx}^2L(\ox,\ov)\ou,\ou\ra=0$, then
\begin{equation*}
\nabla_{xx}^2 L(\ox,\ov)\ou+\nabla\Phi(\ox)^*\bar\eta=0\;\mbox{ for some }\;\bar\eta\in{\cal K}(\oz,\ov)^*\cap\big\{\nabla\Phi(\ox)\ou\big\}^\bot.
\end{equation*}}
To verify this claim, consider the constrained optimization problem:
\begin{equation}\label{ty01}
\mbox{minimize}_{u\in\R^n}\;\frac{1}{2}\la\nabla_{xx}^2L(\ox,\ov)u,u\ra\;\mbox{ subject to }\;\nabla\Phi(\ox)u\in{\cal K}(\oz,\ov).
\end{equation}
It follows from (\ref{re17}) and $\la\nabla_{xx}^2L(\ox,\ov)\ou,\ou\ra=0$ that $\ou$ is an optimal solution to (\ref{ty01}). Using the standard first-order optimality condition and sum rule in \eqref{ty01} yields
\begin{equation*}
0\in\sub_{u}\Big(\Big\la\frac{1}{2}\nabla_{xx}^2L(\ox,\ov)u,u\Big\ra+\dd_{{\cal K}(\oz,\ov)}\big(\nabla\Phi(\ox)u\big)\Big)(\ou)=\nabla_{xx}^2L(\ox,\ov)\ou+\sub_{u}\Big(\dd_{{\cal K}(\oz,\ov)}\big(\nabla\Phi(\ox)u\big)\Big)(\ou).
\end{equation*}
Employing now the calculus rule from Henrion and Outrata \cite[Theorem~5]{ho05} and observing that the calmness assumption therein is automatic due to the linearity of $\nabla\Phi(\ox)u$ and polyhedrality of ${\cal K}(\oz,\ov)$ by Robinson's seminal result from \cite{rob1}, we get
\begin{equation*}
\sub_{u}\Big(\dd_{{\cal K}(\oz,\ov)}\big(\nabla\Phi(\ox)u\big)\Big)(\ou)=\nabla\Phi(\ox)^*N_{{\cal K}(\oz,\ov)}\big(\nabla\Phi(\ox)\ou\big).
\end{equation*}
Substituting it into the above first-order condition gives us the inclusion
\begin{equation*}
0\in\nabla_{xx}^2L(\ox,\ov)\ou+\nabla\Phi(\ox)^*N_{{\cal K}(\oz,\ov)}\big(\nabla\Phi(\ox)\ou\big),
\end{equation*}
which is clearly equivalent to the statement of the claim.\vspace*{0.02in}

Assuming now on the contrary that SOSC in (ii) fails and employing Claim~2, we find the pair $(\xi,\eta):=(\ou,\bar\eta)$ with $\xi\ne0$ satisfying all the conditions in (\ref{crc6}). It says by Theorem~\ref{cric002} that the multiplier $\ov$ is {\em critical} at $\ox$, a contradiction.  This verifies the implication (iii)$\Longrightarrow$(ii).\vspace*{0.03in}

Since the implication (i)$\Longrightarrow$(iv) is trivial, it remains to show that (iv)$\Longrightarrow$(iii) for completing the proof of the theorem. In fact, the equivalence between (iv) and (iii) for $\th\in CPWL$ follows from Theorem~\ref{crit-calm} with $S=S_{\tiny{\rm KKT}}$ and $\Lm=\Lm_{\small{\rm com}}$. We can also verify the implication (i)$\Longrightarrow$(iv) by a direct proof while observing that the isolated calmness of $S_{\tiny{\rm KKT}}$ at $((0,0),(\ox,\ov))$ in (iv) gives us a neighborhood $V$ of $(\ox,\ov)$ such that
\begin{equation*}
S_{\tiny{\rm KKT}}(0,0)\cap V=\big\{(\ox,\ov)\big\}.
\end{equation*}
This implies the existence of a neighborhood $W$ of $\ov$ with $\Lm_{\small{\rm com}}(\ox)\cap W=\{\ov\}$. Since the set $\Lm_{\small{\rm com}}(\ox)$ is convex, we easily deduce from here that $\Lm_{\small{\rm com}}(\ox)=\{\ov\}$, which ensures the validity of \eqref{re15}. Employing finally Theorem~\ref{uplip} in this setting tells us that the unique Lagrange multiplier $\ov$ is {\em noncritical} at $\ox$. This justifies (iii) and thus completes the proof of the theorem. $\h$\vspace*{-0.15in}

\section{Noncriticality, Nondegeneracity, and Robust Isolated Calmness from Lipschitz-Like Property}\sce\vspace*{-0.05in}

The goal of this section is to study relationships between the properties of the KKT solution map $S_{\tiny{\rm KKT}}$ listed in the title and another robust stability property of $S_{\tiny{\rm KKT}}$, which is well-understood and employed in variational analysis and optimization. Recall that a set-valued mapping $F\colon\R^n\tto\R^m$ has the {\em Lipschitz-like/Aubin} (known also as pseudo-Lipschitz) property around $(\ox,\oy)\in\gph F$ if there are neighborhoods $U$ of $\ox$, $V$ of $\oy$ and a number $\ell\ge 0$ such that
\begin{equation}\label{lipl}
F(x_1)\cap V\subset F(x_2)+\ell\|x_1-x_2\|\B\;\mbox{ for all }\;x_1,x_2\in U.
\end{equation}
We know from \cite[Theorem~5.7]{m93} and \cite[Theorem~9.40]{rw} that the latter property can be completely characterized via the following {\em coderivative/Mordukhovich criterion}:
\begin{equation}\label{cod2}
D^*F(\ox,\oy)(0)=\{0\}
\end{equation}
provided that $F$ is closed-graph near $(\ox,\oy)$, where the (limiting) coderivative $D^*$ is defined in \eqref{2.8}. Since the coderivative \eqref{2.8} is {\em robust} and enjoys {\em full calculus}, criterion \eqref{cod2} allows us to efficiently deal with structural mappings that appear in variational analysis and optimization; see, e.g., \cite{m06,rw} and their references for a great many results and applications. We mention a very recent paper \cite{gm16}, where it is shown that the Lipschitz-like property of general constrained systems is implied by another one called the {\em Robinson stability} in \cite{gm16} for which various first-order and second-order sufficient conditions and characterizations are established therein.\vspace*{0.03in}

Let us first deduce from \eqref{cod2} the following description  of the Lipschitz-like property for the solution map $S_{\tiny{\rm KKT}}$ to the KKT system  \eqref{geq4}.\vspace*{-0.05in}

\begin{Proposition}{\bf(equivalent description of the Lipschitz-like property for KKT systems).}\label{morc} Let $(\ox,\ov)\in S_{\tiny{\rm KKT}}(0,0)$ for $S_{\tiny{\rm KKT}}$ from \eqref{s-kkt} with $\th\in CPWL$. Then $S_{\tiny{\rm KKT}}$ is Lipschitz-like around
$((0,0),(\ox,\ov))$ if and only if we have the implication
\begin{equation}\label{lip30}
\left\{\begin{array}{ll}
\nabla_{xx}^2L(\ox,\ov)\xi+\nabla\Phi(\ox)^*\eta=0,\\
\eta\in\big(D^*\sub\th\big)(\oz,\ov)\big(\nabla\Phi(\ox)\xi\big)
\end{array}\right.
\Longrightarrow(\xi,\eta)=(0,0).
\end{equation}
\end{Proposition}\vspace*{-0.05in}
{\bf Proof.} Consider the mapping $G$ from \eqref{F-KKT} with $\Psi=\nabla_x L$. It can be easily checked by the coderivative definition \eqref{2.8} that
$$
(\xi,\eta)\in D^*S_{\tiny{\rm KKT}}\big((0,0),(\ox,\ov)\big)(w_1,w_2)\Longleftrightarrow-(w_1,w_2)\in D^*G\big((\ox,\ov),(0,0)\big)(-\xi,-\eta)
$$
whenever $(\xi,\eta)\in\R^n\times\R^m$ and $(w_1,w_2)\in\R^n\times\R^m$. Employing \cite[Theorem~1.62]{m06} and using the symmetry of the Hessian $\nabla_{xx}^2L(\ox,\ov)$ yield
\begin{equation*}
\begin{array}{lll}
D^*G\big((\ox,\ov),(0,0)\big)(\xi,\eta)&=&\left[\begin{array}{cc}
\nabla_{xx}^2L(\ox,\ov)&-\nabla\Phi(\ox)^*\\\nabla\Phi(\ox)&0\\
\end{array}
\right]\left[\begin{array}{c}
\xi\\\eta\\
\end{array}\right]
+\left[\begin{array}{c}
0\\D^*(\sub\th)^{-1}(\ov,\oz)(\eta)
\end{array}\right]\\\\
&=&\left[\begin{array}{c}
\nabla_{xx}^2L(\ox,\ov)\xi-\nabla\Phi(\ox)^*\eta\\
\nabla\Phi(\ox)\xi+D^*(\sub\th)^{-1}(\ov,\oz)(\eta)
\end{array}
\right].
\end{array}
\end{equation*}
Then \eqref{cod2} tells us that $S_{\tiny{\rm KKT}}$ is Lipschitz-like around $((0,0),(\ox,\ov))$ if and only if
\begin{equation*}
(0,0)\in D^*G\big((\ox,\ov),(0,0)\big)(\xi,\eta)\Longrightarrow(\xi,\eta)=(0,0).
\end{equation*}
Combining this and the above coderivative representation for $G$ ensures description \eqref{lip30}. $\h$\vspace*{0.05in}

We are now in a position to justify that the Lipschitz-like property of $S_{\tiny{\rm KKT}}$ around $((0,0),(\ox,\ov))$ implies that the {\em nondegeneracy  condition} \eqref{fnond} holds. To the best of our knowledge, such a result for multivalued solution maps has been first obtained by Klatte and Kummer \cite[Theorem~1]{kk13} for constrained optimization problems with smooth data. Note that our composite optimization problem \eqref{comp} can be written in the explicit constrained framework \eqref{comp1} but with the {\em nonsmooth} cost. The next theorem derives the nondegeneracy condition \eqref{fnond} from the Lipschitz-like property of \eqref{s-kkt} (and hence the uniqueness of Lagrange multipliers) by a proof different from \cite{kk13} while using some advances of second-order generalized differentiation. Furthermore, in this way we establish {\em noncriticality} of the unique Lagrange multiplier as a consequence of the Lipschitz-like property, which seems to be never mentioned before.\vspace*{-0.1in}

\begin{Theorem}\label{lipnd}{\bf(nondegeneracy and noncriticality from the Lipschitz-like property).}\label{nondeg-ll}
Let $S_{\tiny{\rm KKT}}$ from \eqref{s-kkt} with $\th\in CPWL$ be Lipschitz-like around $((0,0),(\ox,\ov))$. Then we have:

{\bf(i)} $\ox$ is a nondegenerate point of $\Phi$ in the sense of \eqref{fnond}.

{\bf(ii)} There are a neighborhood $O$ of $(0,0)$ and a number $\ve>0$ such that for any $(p_1,p_2)\in O$ the Lagrange multiplier set for the perturbed problem {\rm (\ref{comp3})} defined by
$$
\Lm_{p_1p_2}(x_{p_1p_2}):=\big\{v\in\R^m\big|\;p_1=\nabla_x L(x_{p_1p_2},v),\;v\in\sub\th(\Phi(x_{p_1p_2})+p_2)\big\}
$$
reduces to $\{v_{p_1p_2}\}$, where $(x_{p_1p_2},v_{p_1p_2})\in S_{\tiny{\rm KKT}}(p_1,p_2)\cap \B_\ve(\ox,\ov)$.

{\bf(iii)} $\Lm_{\small{\rm com}}(\ox)=\{\ov\}$, and the multiplier $\ov$ is noncritical.
\end{Theorem}\vspace*{-0.05in}
{\bf Proof.} As discussed in Section~6, the nondegeneracy condition \eqref{fnond} for $\ox$ can be equivalently written as \eqref{nondeg}. To verify the latter, pick $\eta\in\aff\partial\th(\oz)\cap\ker\nabla\Phi(\ox)^*$ and deduce from \cite[Theorem~3.1(ii)]{ms152} that $\eta\in\aff\partial\th(\oz)=(D^* \sub\th)(\oz,\ov)(0)$; thus we come up to
$$
\nabla\Phi(\ox)^*\eta=0\;\mbox{ and }\;\eta\in\big(D^*\sub\th\big)(\oz,\ov)(0).
$$
Since $S_{\tiny{\rm KKT}}$ is Lipschitz-like around $((0,0),(\ox,\ov))$, it follows from Proposition~\ref{morc} that $\eta=0$, and therefore we justify the nondegeneracy assertion (i).

To proceed further with verifying (ii), deduce from \cite[Proposition~4.2]{ms152} that \eqref{fnond} ensures that the set $\Lm_{\small{\rm com}}(\ox)$ is a singleton. Since the Lipschitz-like property is robust/stable under small perturbations of the initial data, we get (ii).

To prove finally (iii), we get from (ii) that $\Lm_{\small{\rm com}}(\ox)=\{\ov\}$, and so it remains to justify the noncriticality of $\ov$. Definition~\ref{crit} requires verifying the implication
\begin{equation*}
0\in\nabla_{xx}^2L(\ox,\ov)\xi+\nabla\Phi(\ox)^*\big(D\sub\th\big)(\oz,\ov)\big(\nabla\Phi(\ox)\xi\big)\Longrightarrow\xi=0.
\end{equation*}
Pick $\xi\in\R^n$ such that $\nabla_{xx}^2L(\ox,\ov)\xi+\nabla\Phi(\ox)^*\eta=0$ for some $\eta\in(D\sub\th)(\oz,\ov)(\nabla\Phi(\ox)\xi)$. Then the derivative-coderivative relationship \eqref{der-cod} yields the conditions
\begin{equation*}
\nabla_{xx}^2L(\ox,\ov)\xi+\nabla\Phi(\ox)^*\eta=0\;\mbox{ and }\;\eta\in\big(D^*\sub\th\big)(\oz,\ov)\big(\nabla\Phi(\ox)\xi\big).
\end{equation*}
The imposed Lipschitz-like property of $S_{\tiny{\rm KKT}}$ around $((0,0),(\ox,\ov))$ tells us by Proposition~\ref{morc} that $\xi=0$, which justifies the noncriticality of $\ov$ and thus completes the proof. $\h$\vspace*{0.05in}

We finish this section by showing that the Lipschitz-like property of $S_{\tiny{\rm KKT}}$ implies the robust isolated calmness of this set. The obtained result can be compared with \cite[Proposition~20 and Corollary~25]{dsz} for problems of constrained optimization with smooth data and nonpolyhedral constraint sets. Recall that our equivalent constrained optimization form \eqref{comp1} of \eqref{comp} intrinsically contains nonsmoothness. The proof presented in \cite{dsz} is based on an involved  result by Fusek \cite{fu} and is different from the second-order variational tools implemented below.\vspace*{-0.1in}

\begin{Theorem}{\bf(robust isolated calmness from the Lipschitz-like property).}\label{lipnd2} If the solution map $S_{\tiny{\rm KKT}}$ from \eqref{s-kkt} with $\th\in CPWL$ enjoys the Lipschitz-like property around $((0,0),(\ox,\ov))$, then it is robustly isolatedly calm at this point.
\end{Theorem}\vspace*{-0.05in}
{\bf Proof.} To justify this result, we argue similarly to the proof of Theorem~\ref{crit-calm} by using now the description of the Lipschitz-like property of  $S_{\tiny{\rm KKT}}$ taken from Proposition~\ref{morc}. Let the pair $(\xi,\eta)\in\R^n\times\R^m$ belong to the set on the left-hand side of \eqref{crc5}. Employing this together with the derivative-coderivative relationship (\ref{der-cod}) says that $(\xi,\eta)$ also belongs to the set on the left-hand side of \eqref{lip30}. Thus the assumed Lipschitz-like property of $S_{\tiny{\rm KKT}}$ tells us that $(\xi,\eta)=(0,0)$ by Proposition~\ref{morc}. Employing finally Theorem~\ref{crit-calm}, we conclude that the solution map $S_{\tiny{\rm KKT}}$ has the isolated calmness property at $((0,0),(\ox,\ov))$, while its robustness is a direct consequence of the implication (iv)$\Longrightarrow$(i) in Theorem~\ref{calm2}.$\h$\vspace*{-0.15in}

\section{Concluding Remarks}\vspace*{-0.05in}

This paper reveals deep connections between {\em critical/noncritical multipliers} for variational systems and {\em second-order generalized differentiation} in variational analysis. We employ second-order constructions in the suggested definition of critical multipliers and then strongly benefit from the recent second-order calculations for the class of extended-real-valued {\em CPWL functions} in terms of their given data. This part exploits the polyhedral epigraphical structure of such functions, which is also used in some proofs based on the Hoffman Lemma. Applications to optimization are done in this paper in the formalism of composite optimization problems that are intrinsically {\em nonsmooth} even if written in the constrained optimization framework.

One of the most important messages for {\em numerical optimization} delivered by obtained results in the class of composite models is that critical multipliers and {\em slow convergence} of major primal-dual algorithms induced by the existence of such multipliers can be {\em ruled out} if we search not arbitrary minimizers but only those satisfying certain {\em stability} properties, which have been recently fully characterized via the problem data. This may allow the user to make some conclusions about algorithm convergence properties {\em a priori} the convergence analysis.

Our future plans concern developing the suggested approach to the study of critical multipliers for variational systems and optimization problems without {\em any polyhedral structure}. Preliminary results confirm the possibility of such developments and their applications to several classes of nonpolyhedral constrained optimization including {\em second-order cone programming}.\\[1ex]
{\bf Acknowledgements.} The first author gratefully acknowledges numerous discussions with Alexey Izmailov and Mikhail Solodov on critical multipliers and related topics. We particularly appreciate sharing with us Izmailov's instructive notes \cite{iz15}. We are also indebted to two anonymous referees and the handling editor for their very careful reading of the paper and making helpful remarks that allowed us to improve the original presentation.
\end{document}